\documentclass[a4paper,fleqn]{cas-dc}

\usepackage[authoryear,longnamesfirst]{natbib}
\usepackage{mathtools}
\newcommand{\redx}[1]{{\color{red}#1}}
\def\tsc#1{\csdef{#1}{\textsc{\lowercase{#1}}\xspace}}
\tsc{WGM}
\tsc{QE}

\newproof{pf}{Proof}
\newtheorem{theorem}{Theorem}
\newtheorem{lemma}{Lemma}%
\newtheorem{proposition}{Proposition}%
\newtheorem{remark}{Remark}%
\usepackage{algorithm, algpseudocode}
\newtheorem{definition}{Definition}
\newtheorem{assumption}{Assumption}
\DeclareMathOperator*{\diag}{diag}

\newcommand*{\minimize}[1]{\operatorname*{minimize}_{#1}\quad} 
 
\newcommand{\subjectto}{\textup{subject to}\quad}

\usepackage[shortcuts]{extdash}
\begin{document}
\let\WriteBookmarks\relax
\def\floatpagepagefraction{1}
\def\textpagefraction{.001}

\shorttitle{Stability-Certified Imitation Learning for Boundary Control of Reaction--Diffusion PDEs}    

\shortauthors{Biazetto et al.}  

\title [mode = title]{Stability-Certified Imitation Learning for Boundary Control of Reaction--Diffusion PDEs}  

\tnotemark[1] 

\tnotetext[1]{This work was partially supported by CAPES under grant 88881.878833/2023-01 (SticAmSud) and CNPq under the grant  310374/2026-0.} 

\author[1]{Paulo Henrique F. Biazetto}[orcid=0000-0002-6497-8687]
\cormark[1]
\ead{paulo.biazetto@posgrad.ufsc.br}


\credit{Methodology, Conceptualization, Software, Writing – Original Draft Preparation}

\affiliation[1]{organization={Department of Automation and Systems Engineering, Federal University of Santa Catarina},
            city={Florianópolis},
            postcode={88040-900}, 
            state={SC},
            country={Brazil}}

\author[2]{Mirko Fiacchini}[orcid=0000-0002-3601-0302]

\ead{mirko.fiacchini@gipsa-lab.fr}


\credit{Methodology, Conceptualization, Supervision, Writing – Review \& Editing}

\affiliation[2]{organization={University Grenoble Alpes, CNRS, Grenoble-INP, GIPSA-Lab},
            city={Grenoble},
            postcode={ 38000}, 
            country={France}}

\author[2]{Christophe Prieur}[orcid=0000-0002-4456-2019]

\ead{christophe.prieur@gipsa-lab.fr}


\credit{Methodology, Conceptualization, Supervision, Writing – Review \& Editing}

\author[1]{Gustavo A. {de Andrade}}[orcid=0000-0001-9560-8833]


\ead{gustavo.artur@ufsc.br}


\credit{Methodology, Conceptualization, Supervision, Writing – Review \& Editing}

\cortext[1]{Corresponding author}



\begin{abstract}
We develop a stability-certified imitation-learning framework for Dirichlet boundary control of one-dimensional reaction-diffusion partial differential equations (PDEs). An eigenfunction expansion of the associated Sturm-Liouville operator decomposes the plant into a finite-dimensional subsystem, containing all unstable and selected stable modes, and an infinite-dimensional residual subsystem. A feedforward neural network is trained to imitate a constrained model predictive control policy on the retained model. A quadratic Lyapunov function and local quadratic constraints on the neural-network activations yield a tractable certificate of regional exponential stability and an ellipsoidal inner approximation of the retained-state region of attraction. To account for spillover, the Lyapunov function is augmented with a spectrally weighted residual-mode term. The resulting composite functional provides sufficient conditions for local exponential stability of the full closed-loop PDE in a negative-order Hilbert space. The finite- and infinite-dimensional requirements are incorporated into an aggregate constrained training formulation and reverified for the final fixed network. Numerical experiments on an unstable reaction-diffusion equation show that the certified neural controller closely reproduces the expert closed-loop behavior while reducing average online evaluation time by approximately three orders of magnitude.
\end{abstract}


\begin{highlights}
\item Stability-certified imitation learning is developed for boundary-controlled PDEs
\item A weighted residual-mode Lyapunov term captures infinite-dimensional spillover
\item LMIs certify regional exponential stability of the complete closed-loop PDE
\item Neural control approximates MPC with substantially lower online evaluation cost
\end{highlights}



\maketitle

\section{Introduction}\label{section:introduction}

Boundary control of partial differential equations (PDEs) must reconcile infinite-dimensional dynamics with stability, performance, and online computational requirements \citep{morris2020control}. Established model-based methods, including backstepping \citep{krstic2008}, semigroup and operator-theoretic techniques \citep{curtain2020}, Lyapunov-based energy methods \citep{bastin2016}, and optimal-control formulations \citep{li1995optimal}, provide systematic tools for controller synthesis and stability analysis. Their guarantees, however, are generally tied to specific analytical representations and controller structures and therefore do not transfer directly to feedback laws represented by general neural networks (NNs). This difficulty is particularly pronounced for boundary-actuated PDEs, for which the boundary input is typically represented by an unbounded control operator in the natural state space.

NN-based methods have been used to identify governing equations and parameters \citep{tanyu2023deep}, approximate PDE solution operators \citep{lu2021deeponet,li2021fourier,kovachki2023neural,luo2025physics}, learn feedback policies through reinforcement learning \citep{mattioni2023enhancing,bhan2024pdecontrolgym}, and approximate control-related operators \citep{bhan2023operator,zhang2026operator}. In particular, neural operators have been employed to approximate the gain-generation and feedback operators arising in PDE backstepping designs, with closed-loop guarantees retained under suitable bounds on the approximation error \citep{bhan2024neuraloperators,krstic2024neuraloperators}. This methodology has also been extended to adaptive boundary control of reaction--diffusion PDEs, where the controller gains are updated as unknown plant parameters are estimated online \citep{bhan2025adaptive}. These results concern learned approximations of prescribed model-based controller structures. A distinct problem arises when the feedback policy itself is learned from demonstrations and represented by a general NN parameterization.

For finite-dimensional systems, stability-constrained imitation learning (IL) combines learning-based control with Lyapunov and robust-control tools. An NN controller can be trained to reduce the discrepancy from an expert policy while enlarging a certified region of attraction (ROA) \citep{yin2021imitation,east2022imitation}. Related formulations consider trade-offs between imitation accuracy and stability margins \citep{pinguet2023data}. These methods commonly employ quadratic constraints (QCs) to characterize the NN activation functions and derive tractable sufficient conditions for closed-loop stability \citep{yin2022stability,zare2024survey}.

Extending this framework to boundary-controlled PDEs requires additional analysis. A stability certificate for a finite-dimensional modal subsystem does not, by itself, establish stability of the complete PDE. Boundary actuation may produce modal input coefficients that are not square summable, and the infinitely many residual modes omitted from controller synthesis may be excited by the learned feedback law. These residual dynamics and the resulting spillover must therefore be included explicitly in the infinite-dimensional stability analysis.

An initial application of stability-constrained IL to NN-based control of a distributed-parameter system was presented in \citep{de2023imitation}. The stability certificate was constructed for a truncated modal representation of the plant, without explicitly including the neglected modes in the Lyapunov functional. Moreover, the synthesis conditions did not impose a prescribed exponential decay rate. These limitations leave open the problem of certifying the complete infinite-dimensional closed-loop system while retaining the flexibility and low online evaluation cost of an imitation-learned policy.

This paper addresses this problem for a class of one-dimensional reaction--diffusion PDEs with Dirichlet boundary actuation. A feedforward NN is trained from state--control demonstrations to reproduce an expert feedback policy. Although the policy is learned from data, the proposed method is not model free: the known PDE dynamics are used to construct and verify the stability certificate. The objective is to obtain a computationally efficient approximation of the expert controller while retaining rigorous closed-loop guarantees.

The main contributions are twofold. First, the PDE is expanded in the eigenfunction basis of the Sturm--Liouville operator associated with the homogeneous boundary conditions, and the resulting modal dynamics are decomposed into a retained subsystem, containing all unstable modes and a finite number of stable modes, and an infinite-dimensional residual subsystem. A quadratic Lyapunov function for the retained dynamics is combined with local QCs describing the NN activation functions and augmented with a spectrally weighted term involving the residual modal coordinates. The weighting accommodates the high-frequency behavior of the modal boundary-input coefficients and defines a composite Lyapunov functional in a negative-order Hilbert space. This construction yields sufficient matrix-inequality conditions for regional exponential stability with a prescribed decay rate, provides a certified ellipsoidal inner approximation of the retained-state ROA, and establishes local exponential stability of the complete infinite-dimensional closed-loop system.

Second, the retained-state stability certificate, state-set containment conditions, and residual-mode spillover requirements are assembled into a single constrained learning formulation. Because the resulting learning problem is nonconvex, convergence of the numerical algorithm does not by itself provide a stability guarantee. After training, the certificate is therefore recomputed using the exact representation of the final fixed NN, and the controller is declared certified only if all finite- and infinite-dimensional verification inequalities are satisfied.

The proposed method is evaluated on an unstable reaction\-/diffusion equation. A constrained model predictive control (MPC) law designed for a sampled-data retained model is used as the expert policy and generates the state--control demonstrations. The learned NN is then implemented as a continuous-time state-feedback controller. The numerical study evaluates the imitation accuracy, the certified retained-state ellipsoid, the prescribed decay rate, the residual-mode response, and the reduction in online evaluation time relative to MPC.

The remainder of the paper is organized as follows. Section~\ref{sec:preliminary} recalls the required results from Sturm--Liouville theory. Section~\ref{sec:prob} introduces the reaction--diffusion system and derives its modal representation. Section~\ref{sec:imitation} presents the stability-certified IL framework. Section~\ref{section:stability_analysis} establishes the infinite-dimensional stability results. Section~\ref{section:ADMM} describes the numerical implementation of the constrained learning procedure. Section~\ref{sec:results} reports the numerical results, and Section~\ref{sec:conclusion}  concludes the paper.

\section{Preliminaries on Sturm--Liouville Theory} \label{sec:preliminary}
Let $L^{2}(0,1)$ be the real Hilbert space of square-integrable functions on the interval $(0,1)$, equipped with the inner product $\langle f,g \rangle_{L^{2}(0,1)} = \int_{0}^{1} f(x)g(x) dx$ and the induced norm $\|f\|_{L^{2}(0,1)}=\sqrt{\langle f,f\rangle_{L^{2}(0,1)}}$.

\begin{definition}\label{def:sturm-liouville}
Let $p\in C^{1}([0,1])$ be a real-valued function satisfying $p(x)>0$, for every $x\in[0,\,1]$. We consider the regular Sturm--Liouville operator $\mathcal{A}: D(\mathcal{A}) \subset L^{2}(0,1) \rightarrow L^{2}(0,1)$ defined by
\begin{equation}
    \mathcal{A}f = -(pf')', \label{eq:SLoperator}
\end{equation}
\noindent for all $f$ in the domain
\begin{align}
    \hspace{-0.8cm} D(\mathcal{A}) &= \{f \in H^{2}(0,1) :\, \cos(\theta_{0})f(0) - \sin(\theta_{0})f'(0) = 0, \nonumber\\
     &~~~~~\cos(\theta_{1})f(1) + \sin(\theta_{1})f'(1) = 0\} \label{eq:bc_1},
\end{align}
\noindent where $\theta_{0},\theta_{1}\in[0,\pi)$.
\end{definition} 

For any $f,g \in D(\mathcal{A})$, Green's identity and the boundary conditions in \eqref{eq:bc_1} give $\langle\mathcal{A}f,g\rangle_{L^{2}(0,1)}-\langle f,\mathcal{A}g\rangle_{L^{2}(0,1)}=\left[p\left(fg'-f'g\right)\right]_{0}^{1}=0$. Thus, $\mathcal{A}$ is symmetric. Standard results for regular Sturm--Liouville operators further imply that $\mathcal{A}$ is self-adjoint, bounded from below, and has a compact resolvent \citep[Theorem~3.2.8]{curtain2020}.

Consequently, the spectrum of $\mathcal{A}$ consists of a countable set of real eigenvalues $\{\lambda_{n}\}_{n=1}^{\infty}$, with no finite accumulation point and with $\lambda_n\to+\infty$ as $n\to\infty$. The corresponding eigenfunctions $\{\Phi_{n}\}_{n=1}^{\infty}$ form a complete orthonormal basis for $L^2(0,1)$ \cite[Theorem 3.2.8]{curtain2020}. Thus, every $f \in L^2(0,1)$ admits the  spectral expansion
\begin{equation}
    f = \sum_{n=1}^{\infty} f_n \Phi_{n}, \quad\text{in }L^2(0,1),
\end{equation}
where $\quad f_n = \langle f, \Phi_{n} \rangle_{L^{2}(0,1)}$.

For every $f \in D(\mathcal{A})$, the action of the operator is given by
\begin{equation}
    \mathcal{A}f = \sum_{n=1}^{\infty} \lambda_{n} f_n \Phi_{n}, \quad\text{in }L^2(0,1).
\end{equation}

Since $\mathcal{A}$ is self-adjoint and bounded from below, $-\mathcal{A}$ generates a ${C}_{0}$-semigroup on $L^{2}(0,1)$. Consequently, for any initial condition $z_{0}\in L^{2}(0,1)$, the homogeneous problem
\begin{align*}
&z_{t}(t,x) = (p(x)z_{x}(t,x))_{x},\\
    &\cos(\theta_{0})z(t,0)-\sin(\theta_{0})z_{x}(t,0)=0,\nonumber\\
    &\cos(\theta_{1})z(t,1)+\sin(\theta_{1})z_{x}(t,1)=0,\nonumber\\
    &z(0,x)=z_{0}(x),\nonumber
\end{align*}
admits a unique mild solution given by
\begin{equation*}
    z(t)=e^{-\mathcal{A}t}z_{0}=\sum_{n=1}^{\infty}e^{-\lambda_{n}t}\langle z_{0},\Phi_{n}\rangle_{L^{2}(0,1)}\Phi_{n},
\end{equation*}
where the series converges in $L^{2}(0,1)$ for every $t\geq 0$. If $z_{0}\in D(\mathcal{A})$, the corresponding solution is a strong solution satisfying $z\in C\big([0,\infty);D(\mathcal A)\big)\cap C^{1}\big([0,\infty);L^2(0,1)\big)$. This spectral representation expresses the PDE dynamics as an infinite-dimensional system of modal ordinary differential equations.

\section{Problem Statement} \label{sec:prob}
We consider the boundary stabilization of an unstable linear reaction--diffusion PDE on a one-dimensional spatial domain. The objective is to construct a computationally efficient NN boundary controller that approximates a constrained expert policy while providing a certificate of local exponential stability for the complete infinite-dimensional closed-loop system. To this end, this section homogenizes the boundary condition and derives the modal representations used for controller synthesis, learning, and stability analysis.

\subsection{The Reaction-Diffusion System}
We consider the following one-dimensional linear reaction\-/diffusion system \citep{lhachemi2020pi}:
\begin{align}
    &z_{t}(t,x) = (p(x)z_{x}(t,x))_{x} + q_c z(t,x), \label{eq:pde_main}\\
    &z_{x}(t,0) = 0, \quad z(t,1) = u(t), \label{eq:pde_bc}\\
    & z(0,x) = z_{0}(x),\label{eq:pde_ini}
\end{align}
\noindent where $t\in[0,\infty)$ is time, $x\in [0,1]$ is space, $p \in {C}^2([0,1])$ with $p(x)>0$, $q_c \in \mathbb{R}$ is the reaction coefficient, $z_{0}\in H^{2}(0,1)$ is the initial condition, and $u:[0,\infty)\rightarrow\mathbb R$ is the boundary control input. We assume that $q_c$ is such that the uncontrolled system has at least one unstable eigenmode. 

The following lemma establishes the well-posedness of \eqref{eq:pde_main}--\eqref{eq:pde_ini} under the assumption that the control input $u\in C^{2}([0,\,\infty))$ \cite{curtain2020}.
\begin{lemma}\label{lemma1}
    Let $p \in C^{1}([0,1])$ with $p(x)>0$ for every $x\in[0,1]$, and let $q_c \in \mathbb{R}$. Additionally, assume that $u\in C^{2}([0,\,\infty))$, and let $z_{0}\in H^{2}(0,1)$ satisfying the compatibility condition
    \begin{align*}
        z_{0}'(0)=0, \qquad z_{0}(1)=u(0).
    \end{align*}

    Then the initial-boundary value problem \eqref{eq:pde_main}--\eqref{eq:pde_ini} admits a unique solution $z\in C([0,\infty);H^2(0,1))\cap C^{1}([0,\infty);\allowbreak L^2(0,1))$.
\end{lemma}

Lemma~\ref{lemma1} concerns classical solutions associated with a prescribed smooth boundary input.  Because the input generated by the NN feedback law does not satisfy this regularity assumption a priori, well-posedness of the feedback interconnection is established separately in the negative-order state space $\mathcal{H}_{-1/2}$, defined in Section~\ref{section:stability_analysis}; see Appendix~\ref{app:well_posedness}.

\subsection{Transformation to a Homogeneous System}

The nonhomogeneous Dirichlet condition \eqref{eq:pde_bc} prevents the homogeneous-domain spectral representation developed in the preceding section from being applied directly to the PDE dynamics. We therefore introduce a change of variables that transforms \eqref{eq:pde_main}--\eqref{eq:pde_ini} into an equivalent system with homogeneous boundary conditions.  For a prescribed input $u\in C^{1}([0,\infty))$, let $b\in C^{2}([0,1])$ satisfy $b'(0)=0$, and $b(1)=-1$, and define
\begin{equation}
w(t,x) = z(t,x) + b(x) u(t).\label{eq:lift}
\end{equation}

Equivalently,
\begin{align}
    z(t,x) = w(t,x) - b(x) u(t). \label{eq:transform}
\end{align}

Substituting \eqref{eq:transform} into \eqref{eq:pde_main}--\eqref{eq:pde_bc} gives
\begin{align}
    w_{t}(t,x) &= (p(x)w_{x}(t,x))_{x} + q_c w(t,x) \nonumber\\
     &\quad + a(x)u(t) + b(x) \dot{u}(t), \label{eq:w_pde} \\
    w_x(t,0) &= 0, \quad w(t,1) = 0, \label{eq:w_bcs}\\
    w(0,x) &= w_{0}(x),
\end{align}
\noindent where $a(x) = - (p(x)b'(x))' - q_c b(x)$ and $w_{0}(x)=z_0(x)+b(x)u(0)$.

Moreover, since $z_0 \in H^2(0,1)$ and $b \in C^2([0,1])$, it follows that $w_{0}\in H^2(0,1)$, given that $u(0)$ is finite. The compatibility conditions $z_{0}'(0)=0$ and $z_{0}(1)=u(0)$, together with the conditions $b'(0)=0$ and $b(1)=-1$, imply that $w_{0}'(0)=0$ and $w_{0}(1)=0$. Hence $w_{0}$ belongs to the domain of the Sturm-Liouville operator $\mathcal{A}$ associated with the homogeneous boundary conditions.

\subsection{Spectral Decomposition}

The homogeneous boundary conditions in \eqref{eq:w_bcs} correspond to a particular case of the Sturm--Liouville operator introduced in Definition~\ref{def:sturm-liouville}. More precisely, by selecting $\theta_0=\frac{\pi}{2}$, $\theta_1=0$, the general boundary conditions in \eqref{eq:bc_1} reduce to $ f'(0)=0$, and $f(1)=0$. Accordingly, the operator associated with the transformed system \eqref{eq:w_pde}--\eqref{eq:w_bcs} is represented by the Sturm--Liouville operator \eqref{eq:SLoperator} with domain
\begin{equation}
    D(\mathcal{A})
    = \left\{  f\in H^2(0,1): f'(0)=0,\quad f(1)=0 \right\}.
    \label{eq:SLdomain}
\end{equation}

Note that for every $f\in D(\mathcal{A})$, integration by parts gives $\langle\mathcal{A}f,f\rangle_{L^{2}(0,1)}=\int_{0}^{1}p(x)|f'(x)|^{2}\,dx$. Since $p(x)>0$ on $[0,1]$ and $f(1)=0$, the operator $\mathcal{A}$ is positive. The spectral properties established in Section~\ref{sec:preliminary} therefore apply, and its eigenvalues can be ordered as 
\begin{align}
    0<\lambda_{1}<\lambda_{2}<\cdots,\qquad\lambda_n\rightarrow+\infty.
\end{align}

The corresponding eigenfunctions $\{\Phi_n\}_{n=1}^{\infty}$ form a complete orthonormal basis for $L^{2}(0,1)$. For every $t\geq 0$, we can project the state $w(t,x)$ onto this eigenbasis in $L^{2}(0,1)$:
\begin{equation*}  
    w(t,x) = \sum_{n=1}^\infty w_n(t)\Phi_{n}(x),
\end{equation*}
where, for every $n\geq 1$ and $t\geq0$,
\begin{align*}
    w_n(t) &= \langle w(t, \cdot), \Phi_{n} \rangle_{L^{2}(0,1)}, \\
    a_{n} &= \langle a, \Phi_{n} \rangle_{L^{2}(0,1)}, \\
    b_{n} &= \langle b, \Phi_{n} \rangle_{L^{2}(0,1)}.
\end{align*}

To find the dynamics of the modal coefficients, we take the time derivative and project onto $\Phi_{n}$:
\begin{align}
    \dot{w}_{n}(t) &= \langle w_{t}(t,\cdot), \Phi_{n} \rangle_{L^{2}(0,1)} \nonumber\\
    &= (-\lambda_{n} + q_{c})w_{n}(t) + a_{n}u(t) + b_{n}\dot{u}(t). \label{eq:w_decomposed}
\end{align}

Let us define the projection coefficients $$z_{n}(t) = \langle z(t,\cdot), \Phi_{n} \rangle_{L^2(0,1)}.$$ 

From the transformation \eqref{eq:transform}, we have 
\begin{equation}
    z_{n}(t) = w_{n}(t) - b_{n}u(t),
\end{equation}
which implies
\begin{align}
    w_{n}(t) = z_{n}(t) + b_{n}u(t). \label{eq:transform_2}
\end{align}

Differentiating \eqref{eq:transform_2} and substituting the result into \eqref{eq:w_decomposed} yields, $\forall n \geq 1$ and $t \geq 0$,
\begin{equation}
    \dot{z}_{n}(t) = (-\lambda_{n}+q_{c}) z_{n}(t) + \beta_{n}u(t), \label{eq:modes}
\end{equation}
\noindent which gives the infinite set of ODEs describing the system dynamics. The coefficients $\beta_{n} = (-\lambda_{n} + q_{c})b_{n} + a_{n}$ characterize how the control input influences each mode.

\begin{remark}
For initial conditions in $\mathcal H_{-1/2}$, the preceding integration-by-parts calculation is understood in the variational sense. More precisely, the nonhomogeneous boundary condition is first lifted, the resulting homogeneous problem is interpreted as a mild evolution equation, and the modal identity is obtained by testing that equation against $\Phi_n$. This argument is detailed in Appendix~A and follows the framework of \cite[Remark~4 and Appendix]{Lhachemi26}.
\end{remark}

By Parseval's identity, $\|z(t,\cdot)\|_{L^2(0,1)}^2 = \sum_{n=1}^{\infty}|z_n(t)|^2$. Consequently, $z(t,\cdot)\in L^2(0,1)$ if and only if the sequence of modal coefficients $\{z_n(t)\}_{n\geq1}$ belongs to $\ell^2$. Thus, the space $\ell^2$ provides an equivalent modal representation of the original state space $L^2(0,1)$.

The coefficients $\beta_n$ in \eqref{eq:modes} characterize the influence of the boundary control input on each eigenmode. The following lemma establishes their boundary representation and the properties needed in the subsequent analysis.

\begin{lemma} \label{lemma:beta_properties}
The control effectiveness coefficients $\{\beta_{n}\}_{n=1}^{\infty}$ possess the following three fundamental properties:
\begin{enumerate}
    \item \label{lemma:boundary} They admit the boundary representation
    \begin{equation}
        \beta_n = p(1)\left[b(1)\Phi_{n}'(1) - b'(1)\Phi_{n}(1)\right]. \label{eq:beta_formula}
    \end{equation}

    \item \label{lemma:controllability} The control input  acts nontrivially on every mode, i.e.,
    \begin{equation}
        \beta_n \neq 0, \quad \forall n \ge 1. \label{eq:beta_nonzero}
    \end{equation}

    \item \label{lemma:asymtotic} The coefficients satisfy the high-frequency estimate
    \begin{equation}
        \left|\beta_{n}\right| = O(\sqrt{\lambda_{n}}) \quad \text{as} \quad n \to \infty. \label{eq:beta_asymptotic}
    \end{equation}
\end{enumerate}
\end{lemma}

\begin{pf}
We prove each property in turn.

\textit{\ref{lemma:boundary} Boundary representation \eqref{eq:beta_formula}:}
From the definitions of $a_n$ and $a$, we have
\begin{equation*}
    a_n=\left\langle -(pb')'-q_cb,\Phi_n\right\rangle_{L^2(0,1)}.
\end{equation*}

Applying Green's identity yields
\begin{align*}
    a_n&=\left[ p\left(b\Phi_n'-b'\Phi_n\right)\right]_{0}^{1}+\left\langle b,\mathcal{A}\Phi_n \right\rangle_{L^2(0,1)} -q_cb_n.
\end{align*}

Since $b'(0)=0$, $\Phi_n'(0)=0$, and $\mathcal{A}\Phi_n=\lambda_n\Phi_n$, it follows that
\begin{equation*}
    a_n=p(1)\left[b(1)\Phi_n'(1)-b'(1)\Phi_n(1)\right]+(\lambda_n-q_c)b_n.
\end{equation*}

Substituting this expression into $\beta_n=(-\lambda_n+q_c)b_n+a_n$ gives \eqref{eq:beta_formula}.

\textit{\ref{lemma:controllability} Nonvanishing of the modal coefficients \eqref{eq:beta_nonzero}}:
We proceed by contradiction. Suppose that $\beta_n=0$ for some $n\geq1$. From \eqref{eq:beta_formula}, we have $b(1)\Phi_n'(1)-b'(1)\Phi_n(1)=0$, since $p(1)>0$.

Because $\Phi_n\in D(\mathcal{A})$, the boundary condition at $x=1$ gives $\Phi_n(1)=0$. Hence, $b(1)\Phi_n'(1)=0$. Since the lifting function satisfies $b(1)=-1$, it follows that $\Phi_n'(1)=0$. Therefore, $\Phi_n(1)=0$, and $\Phi_n'(1)=0$.

The eigenfunction $\Phi_n$ satisfies the second-order differential equation
\begin{equation*}
    -(p(x)\Phi_n'(x))' =\lambda_n\Phi_n(x).
\end{equation*}

Since $p\in C^{1}([0,1])$ and $p(x)>0$, the corresponding initial-value problem with data prescribed at $x=1$ has a unique solution. The zero data above therefore imply $\Phi_n\equiv0$, contradicting the fact that $\Phi_n$ is an eigenfunction. Consequently, $ \beta_n\neq0$, $\forall n\geq 1$.

\textit{\ref{lemma:asymtotic} Asymptotic behavior \eqref{eq:beta_asymptotic}}:
Since $\Phi_n(1)=0$ and $b(1)=-1$, the boundary representation \eqref{eq:beta_formula} reduces to
\begin{equation*}
    \beta_n=-p(1)\Phi_n'(1).
\end{equation*}

For the normalized eigenfunctions of the regular Sturm--Liouville problem, the endpoint derivatives satisfy \cite[Section~8.2]{tucsnak2009observation}
\begin{equation*}
    |\Phi_n'(1)|= O\left(\sqrt{\lambda_n}\right),\qquad\text{as }n\to\infty.
\end{equation*}
Therefore,
\begin{equation*}
    |\beta_n|=p(1)|\Phi_n'(1)|=O\left(\sqrt{\lambda_n}\right),
\end{equation*}
which proves \eqref{eq:beta_asymptotic}. \hfill $\square$
\end{pf}

Item~\ref{lemma:boundary} of Lemma~\ref{lemma:beta_properties} shows that the integral expression defining $\beta_n$ can be reduced to a boundary term evaluated at $x=1$, where the control input is applied. Moreover, Item~\ref{lemma:controllability} establishes that $\beta_n\neq0$ for every $n\geq1$. Hence, the boundary input acts
nontrivially on every eigenmode and, in particular, on all potentially unstable modes.

The estimate in Item~\ref{lemma:asymtotic} shows that the modal control coefficients grow at most proportionally to $\sqrt{\lambda_n}$. Although this sequence need not be square-summable, the high-frequency dissipation of the parabolic system regularizes its effect on the residual modes. The following lemma shows that the scaled coefficient sequence $\left\{\frac{\beta_n}{\lambda_n}\right\}_{n\geq1}$ belongs to $\ell^2$ and quantifies the corresponding residual-mode contribution.

\begin{lemma}
\label{lem:effective_spillover}
Consider the modal dynamics \eqref{eq:modes}.  Assume that $z(0,\cdot)\in L^2(0,1)$, $u\in L^\infty(0,\infty)$  and that there exist constants $\delta_0 > 0$ and $n_s\in\mathbb N$ such that
\begin{equation}
    -\lambda_n+q_c\leq-\delta_0<0,\qquad \forall n\geq n_s+1. \label{eq:tail_stability}
\end{equation}

Assume further that, for all sufficiently large $n$, $|\beta_n|\leq c_\beta\lambda_n^\kappa$ for some $c_\beta>0$ and $\kappa\leq\tfrac12$, and that $\frac{\lambda_n}{n^2}\longrightarrow c_\lambda>0$ as $n\to\infty$.

Then there exists $n_s^\ast\geq n_s$ such that, for every $n_0\geq n_s^\ast$, the residual state $z_{\mathrm{tail}}(t)\coloneqq\{z_n(t)\}_{n=n_0+1}^{\infty}$ belongs to $\ell^2$ for every $t\geq0$ and satisfies
\begin{align}
    \|z_{\mathrm{tail}}(t)\|_{\ell^2}^2 \leq{}&2\sum_{n=n_0+1}^{\infty}e^{-2(\lambda_n-q_c)t}z_n(0)^2
    \nonumber\\
    &+8\left( \sum_{n=n_0+1}^{\infty}\frac{\beta_n^2}{\lambda_n^2}\right)\|u\|_{L^\infty(0,\infty)}^2.
    \label{eq:effective_spillover}
\end{align}

Moreover, $S_2(n_0)\coloneqq\sum_{n=n_0+1}^{\infty}\frac{\beta_n^2}{\lambda_n^2} <\infty$, $S_2(n_0)\longrightarrow 0$ as $n_0\to\infty$.
\end{lemma}
\begin{pf}
For every $n\geq n_0+1$, the variation-of-constants formula gives
\begin{align}
    \hspace{-0.75cm} z_n(t) = e^{-(\lambda_n-q_c)t}z_n(0)+\beta_n\int_0^t e^{-(\lambda_n-q_c)(t-s)}u(s)\,ds. \label{eq:solution_tail}
\end{align}

Since $\lambda_n\to+\infty$, there exists $n_s^\ast\geq n_s$ such that $\lambda_n-q_c\geq\lambda_n/2$ for all $n\geq n_s^\ast+1$. Hence, for $n_0\geq n_s^\ast$,
\begin{align*}
    |z_n(t)|\leq e^{-(\lambda_n-q_c)t}|z_n(0)|+\frac{2|\beta_n|}{\lambda_n}\|u\|_{L^\infty(0,\infty)}.
\end{align*}

Using the fact that for any $a,b\in\mathbb{R}$, $(a+b)^2\leq2a^2+2b^2$ and summing over $n\geq n_0+1$ yields
\begin{align*}
    \|z_{\mathrm{tail}}(t)\|_{\ell^2}^2 \leq{}&
    2\sum_{n=n_0+1}^{\infty}
    e^{-2(\lambda_n-q_c)t}|z_n(0)|^2\\
    &+8\|u\|_{L^\infty(0,\infty)}^2\sum_{n=n_0+1}^{\infty}\frac{\beta_n^2}{\lambda_n^2}.
\end{align*}

Finally, since $\lambda_n/n^2\to c_\lambda>0$, there exist $\underline c_\lambda>0$ and $N_\lambda$ such that $\lambda_n\geq\underline c_\lambda n^2$ for all $n\geq N_\lambda$. Therefore, if $|\beta_n|\leq c_\beta\lambda_n^\kappa$, it follows that $\frac{\beta_n^2}{\lambda_n^2} \leq Cn^{4\kappa-4}=O(n^{-2}),$ because $\kappa\leq\tfrac12$. Therefore, $\sum_{n=1}^{\infty}\frac{\beta_n^2}{\lambda_n^2}<\infty$, and the corresponding series tails converge to zero as $n_0\to\infty$, proving that $z_{\mathrm{tail}}(t)\in\ell^2$ for every $t\geq0$.

This concludes the proof. \hfill $\square$
\end{pf}

\subsection{Finite-Dimensional Model}
For controller design, we retain a finite number of modes of the infinite-dimensional system \eqref{eq:modes}. Let $n_0\ge n_s$, where $n_s$ is defined in Lemma~\ref{lem:effective_spillover}. Then all modes with $n\ge n_0+1$ satisfy $-\lambda_n+q_c\le -\delta_0<0$, and are exponentially stable. Hence, the retained subsystem contains all unstable modes and may also include slowly decaying stable modes. The residual modes are excluded from the controller parameterization but remain explicitly included in the infinite-dimensional stability analysis.

The truncated finite-dimensional model is
\begin{align}
    \dot{Z}_{n_{0}}(t) &= A_{n_{0}}Z_{n_{0}}(t) + B_{n_{0}}u(t), \label{eq:sys_trunc}
\end{align}
\noindent where 
\begin{align*}
    \hspace{-0.75cm}Z_{n_{0}}(t) &= \begin{bmatrix} z_{1}(t) & z_{2}(t) & \cdots & z_{n_{0}}(t) \end{bmatrix}^{\!\top}, \\
    A_{n_{0}} &= \diag\left(-\lambda_{1}+q_{c},\, -\lambda_{2}+q_{c},\, \dots ,\, - \lambda_{n_{0}}+q_{c} \right), \\
    B_{n_{0}} &= \begin{bmatrix} \beta_{1} & \beta_{2} & \cdots & \beta_{n_{0}} \end{bmatrix}^{\!\top}.
\end{align*}

The corresponding initial condition is $$ Z_{n_{0}}(0) = Z_0 = (\langle z_0,\Phi_1\rangle,\dots,\langle z_0,\Phi_{n_0}\rangle)^{\top}.$$

\section{Imitation Learning} \label{sec:imitation}

This section introduces the NN policy and the algebraic representations required for the stability-constrained IL formulation. In particular, the NN nonlinearities are isolated and characterized through QCs for use in the subsequent Lyapunov analysis and constrained training procedure.

Throughout this section, $|\cdot|$ denotes the Euclidean norm for vectors and its induced norm for matrices.

\subsection{Neural Network Controller}
The controller is parameterized by a static feedforward NN policy $\pi:\mathbb{R}^{n_0}\rightarrow\mathbb{R}$ with $\ell\geq 1$ hidden layers. For an input $Z_{n_0}\in\mathbb{R}^{n_0}$, define the layer maps recursively as
\begin{align}
    \hspace{-0.75cm}\omega_0(Z_{n_0}) &= Z_{n_0},\label{eq:nn_controller1}\\
    \hspace{-0.75cm}\omega_i(Z_{n_0}) &=\phi_i\bigl(W_i\omega_{i-1}(Z_{n_0})+b_i\bigr),\quad i\in\{1,\ldots,\ell\}, \label{eq:nn_controller2}\\
    \hspace{-0.75cm}\pi(Z_{n_0}) &= W_{\ell+1}\omega_\ell(Z_{n_0})+b_{\ell+1}.
    \label{eq:nn_controller3}
\end{align}

Let \(n_i\geq1\) denote the width of the \(i\)-th hidden layer. All vectors are understood to be column vectors. The network parameters satisfy
\[
    W_i\in\mathbb{R}^{n_i\times n_{i-1}},
    \qquad
    b_i\in\mathbb{R}^{n_i},
    \qquad i\in\{1,\ldots,\ell\},
\]
and
\[
    W_{\ell+1}\in\mathbb{R}^{1\times n_\ell},
    \qquad
    b_{\ell+1}\in\mathbb{R}.
\]

Let \(\varphi:\mathbb{R}\to\mathbb{R}\) denote the scalar activation function. Its componentwise extension to the \(i\)-th hidden layer is
\[
    \phi_i:\mathbb{R}^{n_i}\to\mathbb{R}^{n_i},
    \qquad
    \phi_i(\nu_i)
    \coloneqq
    \begin{bmatrix}
    \varphi(\nu_{i,1})\\
    \vdots\\
    \varphi(\nu_{i,n_i})
    \end{bmatrix}.
\]

Thus, the recursively defined layer maps satisfy
\[
    \omega_i:\mathbb{R}^{n_0}\to\mathbb{R}^{n_i},
    \qquad i\in\{1,\ldots,\ell\}.
\]

The same scalar activation function \(\varphi\), chosen here as the hyperbolic tangent, is used in every hidden layer.

The feedback input is obtained by evaluating the policy along the retained-state trajectory:
\begin{equation*}
    u(t)=\pi\bigl(Z_{n_0}(t)\bigr).
\end{equation*}

Similarly, the time-dependent signal generated by the $i$-th layer is
\begin{equation*}
    \omega_i(t)\coloneqq \omega_i\bigl(Z_{n_0}(t)\bigr),\qquad i\in\{0,\ldots,\ell\}.
\end{equation*}

Thus, $\omega_i(Z_{n_0})$ denotes the layer map evaluated at a generic input, whereas $\omega_i(t)$ denotes its value along the closed-loop trajectory.

To preserve the origin as an equilibrium, impose
\begin{equation}
    \pi(0)=0.
    \label{eq:pi_origin}
\end{equation}

Equivalently, parameterize the output bias as
\begin{equation*}
    b_{\ell+1} =- W_{\ell+1}\omega_\ell(0_{n_0}),
\end{equation*}
where $0_{n_0}$ denotes the zero vector in $\mathbb{R}^{n_0}$. Here, $\omega_\ell(0_{n_0})$ means evaluation at the zero network input, not evaluation at time $t=0$.

\begin{lemma}\label{lemma:lipschitz_nn}
Suppose that, in \eqref{eq:nn_controller1}--\eqref{eq:nn_controller3}, the activation function at each hidden layer is the componentwise hyperbolic tangent. Then the policy $\pi:\mathbb R^{n_0}\to\mathbb R$ defined by \eqref{eq:nn_controller1}--\eqref{eq:nn_controller3} is of class $C^\infty$, globally Lipschitz, and globally bounded. More precisely,
\begin{align*}
    |\pi(Z)-\pi(\hat Z)| &\leq L_\pi\|Z-\hat Z\|, \qquad \forall Z,\hat Z\in\mathbb R^{n_0}, \\
    |\pi(Z)|&\leq M_\pi,\qquad\forall Z\in\mathbb R^{n_0},
\end{align*}
where one may take
\begin{equation*}
    L_\pi=\|W_{\ell+1}\|\prod_{i=1}^{\ell}\|W_i\|,\quad M_\pi = \|W_{\ell+1}\|\sqrt{n_\ell} + |b_{\ell+1}|.
\end{equation*}
\end{lemma}

\begin{pf}
Each affine map $\omega_{i-1}\mapsto W_i \omega_{i-1} +b_i$ is of class $C^\infty$, and the scalar function $\tanh:\mathbb{R}\to\mathbb{R}$ is also of class $C^\infty$. Since finite compositions of $C^\infty$ maps remain $C^\infty$, it follows that $\pi\in C^\infty(\mathbb{R}^{n_0};\mathbb{R})$.

Moreover, since $\frac{d}{ds}\tanh(s)=\operatorname{sech}^2(s)$ and $0\leq\operatorname{sech}^2(s)\allowbreak\leq1$ for every $s\in\mathbb R$, the mean value theorem implies that $\tanh$ is globally $1$-Lipschitz. Hence, for every
$i=1,\dots,\ell$,
\begin{align*}
    \|\omega_i(Z)-\omega_i(\hat Z)\| \le \|W_i\|\,\|\omega_{i-1}(Z)-\omega_{i-1}(\hat Z)\|.
\end{align*}

Iterating this estimate along the layers yields
\begin{align*}
    \|\omega_\ell(Z)-\omega_\ell(\hat Z)\| \le \Bigl(\prod_{i=1}^{\ell}\|W_i\|\Bigr)\|Z-\hat Z\|.    
\end{align*}

Therefore,
\begin{align*}
    |\pi(Z)-\pi(\hat Z)|
    &\le \|W_{\ell+1}\|\,\|\omega_\ell(Z)-\omega_\ell(\hat Z)\|\\
    &\le \|W_{\ell+1}\|\Bigl(\prod_{i=1}^{\ell}\|W_i\|\Bigr)\|Z-\hat Z\|.
\end{align*}

Thus $\pi$ is globally Lipschitz with 
\begin{equation*}
    L_\pi=\|W_{\ell+1}\|\prod_{i=1}^{\ell}\|W_i\|.
\end{equation*} 

Moreover, since the hyperbolic tangent takes values in $(-1,1)$, one
has
\begin{equation*}
    \|\omega_\ell(Z)\|\leq\sqrt{n_\ell}, \qquad \forall Z\in\mathbb R^{n_0}.
\end{equation*}

Therefore,
\begin{align*}
    |\pi(Z)| &\leq \|W_{\ell+1}\|\,\|\omega_\ell(Z)\| +|b_{\ell+1}| \nonumber\\
    &\leq \|W_{\ell+1}\|\sqrt{n_\ell} + |b_{\ell+1}|.
\end{align*}

Thus, with $M_\pi =\|W_{\ell+1}\|\sqrt{n_\ell} + |b_{\ell+1}|$ the policy $\pi$ is globally bounded. \hfill $\square$
\end{pf}

\subsection{Isolation of NN Nonlinearities}

To facilitate the Lyapunov and QC analysis, the NN controller is rewritten as the feedback interconnection of a linear system and a static nonlinearity. This representation isolates the activation functions from the affine layers and follows the standard framework used in robust control.

All layer quantities in this subsection are evaluated along the retained-state trajectory. In particular, $\omega_i(t)=\omega_i\bigl(Z_{n_0}(t)\bigr)$.

For $i\in\{1,\dots,\ell\}$, let $\nu_{i}\in\mathbb{R}^{n_{i}}$ be the input to the activation function $\phi_{i}$, that is,
\begin{align*}
    \nu_{i}(t) = W_{i}\omega_{i-1}(t) + b_{i},\qquad i\in\{1,\dots,\ell\}.
\end{align*}

Then, the operation of the $i$-th layer can be expressed as $\omega_{i}(t) = \phi_{i}(\nu_{i}(t))$. Stacking the inputs and outputs of all activation functions, we get
\begin{align*}
   \hspace{-0.5cm} \nu_{\phi}(t) \triangleq \begin{bmatrix}
    \nu_{1}(t)\\
    \vdots\\
    \nu_{\ell}(t)
    \end{bmatrix} \in \mathbb{R}^{n_{\phi}}, \quad \mbox{and} \quad \omega_{\phi} (t)\triangleq \begin{bmatrix}
    \omega_{1}(t)\\
    \vdots\\
    \omega_{\ell}(t)
    \end{bmatrix} \in \mathbb{R}^{n_{\phi}},
\end{align*}
\noindent where $n_{\phi}\triangleq \sum_{i=1}^{\ell}n_{i}$. 

Now, define the vector of stacked activation functions $\phi:\mathbb{R}^{n_{\phi}}\rightarrow \mathbb{R}^{n_{\phi}}$ such that
\begin{align}
    \phi(\nu_{\phi}) \triangleq \begin{bmatrix}
    \phi_{1}(\nu_{1})\\
    \vdots\\
    \phi_{\ell}(\nu_{\ell})
    \end{bmatrix}. \label{eq:vector_nonlinearities}
\end{align}

In this framework, the NN control policy $\pi$ can be represented as the following
augmented linear interconnection:
\begin{align}
    \begin{bmatrix}
    u(t) \\ \nu_{\phi}(t)
    \end{bmatrix} &= \mathcal{N} \begin{bmatrix}
    Z_{n_{0}}(t) \\ \omega_{\phi}(t) \\ 1
    \end{bmatrix}, \label{eq:nn_framework1} \\
    \omega_{\phi}(t) &= \phi(\nu_{\phi}(t)), \label{eq:nn_framework2}
\end{align}
\noindent where
\begin{align*}
    \mathcal{N} &= \left[\begin{array}{c|cccc|c}
    0 & 0 & 0 & \hdots & W_{\ell+1} & b_{\ell+1}\\\hline
    W_{1} & 0 & \hdots & 0 & 0 & b_{1}\\
    0 & W_{2} & \hdots & 0 & 0 & b_{2}\\
    \vdots & \vdots & \ddots & \vdots & \vdots & \vdots\\
    0 & 0 & \hdots & W_{\ell} & 0 & b_{\ell}
    \end{array}\right] \\
    &= \begin{bmatrix}
    N_{u z} & N_{u\omega} & N_{ub}\\
    N_{\nu z} & N_{\nu \omega} & N_{\nu b}
    \end{bmatrix}.
\end{align*}

The decomposition \eqref{eq:nn_framework1}--\eqref{eq:nn_framework2} isolates the activation-function nonlinearities, allowing their input--output behavior to be described by QCs and incorporated into the subsequent stability analysis.

\subsection{Quadratic Constraints}
To analyze the stability of the NN-controlled system, the activation nonlinearities are characterized through QCs. The resulting QC representation replaces the nonlinear activation functions by algebraic inequalities relating the pre- and post-activation signals, allowing the closed-loop system to be analyzed using Lyapunov techniques.

In the scalar case, we say that a nonlinear function $\varphi: \mathbb{R}\rightarrow\mathbb{R}$, with $\varphi(0) = 0$, is sector bounded in the sector $[m,\,r]$, with $m \leq r < \infty$, if the following condition holds for all $\nu \in \mathbb{R}$:
\begin{align}
    (\varphi (\nu ) - m \nu ) ( r \nu - \varphi (\nu ) ) \geq 0.\label{eq:local_sc}
\end{align}



Since the activation functions act independently on each neuron, the scalar sector conditions can be stacked to obtain a quadratic constraint for the vector-valued mapping $\phi$. Indeed, for each \(k\in\{1,\ldots,n_\phi\}\), let \(\nu_{\phi,k}\) and \(\omega_{\phi,k}=\varphi(\nu_{\phi,k})\) denote the input and output of the \(k\)-th scalar activation. 

    \redx{Let \(m_\phi \coloneqq (m_1,\ldots,m_{n_\phi})^\top\) and \(r_\phi \coloneqq (r_1,\ldots,r_{n_\phi})^\top\) collect the lower and upper sector bounds, respectively, with \(m_k \leq r_k\) for each \(k\). Define the corresponding diagonal sector matrices by
    \begin{equation*}
        M_\phi \coloneqq \operatorname{diag}(m_\phi), \qquad
        R_\phi \coloneqq \operatorname{diag}(r_\phi).
    \end{equation*}
    }

The stacked activation mapping is said to satisfy the componentwise sector bounds \([m_\phi,r_\phi]\) if
\begin{align}
    (\omega_{\phi,k}-m_k\nu_{\phi,k})
    (r_k\nu_{\phi,k}-\omega_{\phi,k})\geq0
    \label{eq:local_sc_matrix}
\end{align}
for every \(k\) and every admissible \(\nu_\phi\).

Expanding \eqref{eq:local_sc_matrix} yields the equivalent quadratic constraint
\begin{align}
    \begin{bmatrix} \nu_{\phi,k} \\ \omega_{\phi,k} \end{bmatrix}^{\!\top}
    \begin{bmatrix} -2m_kr_k&m_k+r_k \\ m_k+r_k&-2 \end{bmatrix}
    \begin{bmatrix} \nu_{\phi,k} \\ \omega_{\phi,k} \end{bmatrix}
    \geq0. \label{eq:ineq_sector_condition}
\end{align}

For the NN controller \eqref{eq:nn_framework1}--\eqref{eq:nn_framework2} the local sector conditions can be obtained by a numerical procedure given the set of weight matrices $W_i$, the bias vectors $b_{i}$, and monotonically non-decreasing activation functions, $\phi_{i}$. To compute local sector bounds, we introduce the admissible state set
\begin{align}
    \mathcal Z = \{Z_{n_0}\in\mathbb R^{n_0}: -s\le SZ_{n_0}\le s\},\label{eq:admissible_set}
\end{align}
where $S\in\mathbb R^{n_z\times n_0}$ and $s\in\mathbb R_{>0}^{n_z}$.

The corresponding bounds on the pre-activation variables are obtained by propagating interval bounds through the affine layers of the neural network. Since the activation function $\tanh$ is monotone, these bounds can be computed recursively from the first to the last hidden layer. The resulting local sector coefficients are then determined using the procedure proposed in \cite{yin2021imitation, de2023imitation}, which is summarized in the following lemma.

\begin{lemma} \label{lemma:sector_condition}
    Let $\mathcal{Z}_{h} = \{Z_{n_{0}} \in\mathbb{R}^{n_{0}}: \underline{Z}\leq Z_{n_{0}} \leq \overline{Z}\}$ be the smallest hypercube such that $\mathcal{Z} \subseteq \mathcal{Z}_{h}$. Additionally, consider that the activation functions of the NN defined in \eqref{eq:nn_framework1}--\eqref{eq:nn_framework2} are given by the hyperbolic tangent function. 
    Then, there exist vectors \(m_\phi,r_\phi\in\mathbb R^{n_\phi}\) such that, for every \(Z_{n_0}\in\mathcal Z_h\), each scalar activation input-output pair satisfies its corresponding local sector condition \eqref{eq:local_sc_matrix}.
\end{lemma}

Consequently, whenever \(Z_{n_0}(t)\in\mathcal Z_h\), all componentwise quadratic constraints \eqref{eq:ineq_sector_condition} hold simultaneously.


\subsection{Equivalent Closed-Loop Representation} 

The componentwise QC framework developed in the preceding subsection provides a quadratic description of each activation input-output pair. However, when the NN controller is synthesized jointly with a Lyapunov certificate, the resulting stability conditions are nonconvex because the admissible preactivation intervals, and hence the componentwise local sector matrices $M_\phi=M_\phi(\mathcal N)$ and $R_\phi=R_\phi(\mathcal N)$, depend on the NN parameters.

\begin{figure}
    \centering
    \includegraphics[width=0.9\linewidth]{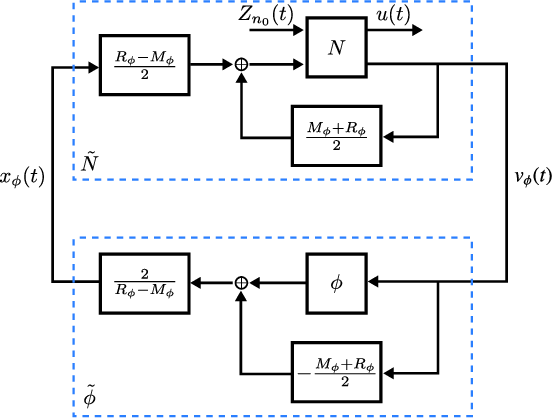}
    \caption{Normalization loop.}
    \label{fig:NN-regularization}
\end{figure}

For fixed local sector matrices $M_\phi$ and $R_\phi$, we employ the loop transformation shown in Figure~\ref{fig:NN-regularization}. This transformation embeds the local sector information into an equivalent linear interconnection while normalizing each activation nonlinearity to the fixed componentwise sector $[-1,1]$. This representation separates the fixed sector information from the normalized nonlinearities and facilitates the formulation of the subsequent stability conditions as LMIs using nonnegative diagonal multipliers.

For notational simplicity, the explicit time dependence of the signals is omitted whenever no ambiguity arises. The following lemma formalizes this equivalent representation.

\begin{lemma}\label{lemma:loop_transformation}
Consider the NN policy $\pi$ represented by \eqref{eq:nn_framework1}--\eqref{eq:nn_framework2} in feedback with the finite-dimensional plant \eqref{eq:sys_trunc}, and assume that the closed-loop system has an equilibrium point $(Z_{n_0}^*,u^*,\nu_\phi^*,\omega_\phi^*)$. After shifting this equilibrium to the origin, 
assume that the activation mapping satisfies the componentwise sector conditions \eqref{eq:local_sc_matrix} over the admissible set \(\mathcal Z\). The translated NN is characterized by the matrix
\begin{equation}
    N \coloneqq
    \begin{bmatrix}
        N_{uz} & N_{u\omega}\\
        N_{\nu z} & N_{\nu\omega}
    \end{bmatrix},
    \label{eq:nn_translated}
\end{equation}
that is,
\begin{equation*}
    \begin{bmatrix}
        u\\
        \nu_\phi
    \end{bmatrix}
    = N
    \begin{bmatrix}
        Z_{n_0}\\
        \omega_\phi
    \end{bmatrix}.
\end{equation*}

Define the center and radius matrices 
\begin{equation*} 
    c_\phi \coloneqq \frac{M_{\phi}+R_{\phi}}{2}, \qquad \rho_\phi \coloneqq \frac{R_{\phi}-M_{\phi}}{2}, 
\end{equation*} 
\noindent and the intermediate block matrices $C_1 \coloneqq N_{u\omega}\rho_\phi$,  $C_2 \coloneqq N_{u\omega}c_\phi$, $C_3 \coloneqq N_{\nu\omega}\rho_\phi$, and  $C_4 \coloneqq N_{\nu\omega}c_\phi$. 

Assume that $R_\phi-M_\phi\succ0$, so that $\rho_\phi=(R_\phi-M_\phi)/2$ is invertible. Then, under the loop transformation 
\begin{equation} 
    \omega_{\phi}(t) = \rho_\phi\, x_{\phi}(t) + c_\phi\, \nu_{\phi}(t), \label{eq:loop_transform} 
\end{equation} 
\noindent the translated NN controller admits an equivalent controller representation defined by a matrix $\tilde{N}\in\mathbb{R}^{(1+n_\phi)\times(n_0+n_\phi)}$ and a normalized nonlinearity $\tilde\phi$:
\begin{align} 
    \begin{bmatrix} 
    u \\ \nu_{\phi}
    \end{bmatrix} &= \tilde{N} 
    \begin{bmatrix} Z_{n_{0}} \\ 
    x_{\phi}
    \end{bmatrix}, \label{eq:nn_tilde_framework1} \\ 
    x_{\phi}(t) &= \tilde{\phi}(\nu_{\phi}(t)), \label{eq:nn_tilde_framework2} 
\end{align} 
\noindent where 
\begin{align} 
    \tilde{\phi}(\nu_{\phi}(t)) &\coloneqq \rho_\phi^{-1}\big(\phi(\nu_{\phi}(t))-c_\phi\nu_{\phi}(t)\big),
    \label{eq:tildephi_def}\\
    \tilde{N} &\coloneqq \begin{bmatrix} \tilde{N}_{u z} & \tilde{N}_{u x} \\ 
    \tilde{N}_{\nu z} & \tilde{N}_{\nu x} 
    \end{bmatrix}, \label{eq:nn_normalized} 
\end{align} 
\noindent with
\begin{align*}
    \tilde{N}_{u z} &= N_{uz} + C_{2}(I_{n_\phi} - C_{4})^{-1}N_{\nu z},\\
    \tilde{N}_{u x} &=  C_{1}+C_{2}(I_{n_\phi}-C_{4})^{-1}C_{3},\\
    \tilde{N}_{\nu z} &= (I_{n_\phi}-C_{4})^{-1}N_{\nu z},\\
    \tilde N_{\nu x} &= (I_{n_\phi}-C_{4})^{-1}C_{3}. 
\end{align*}

The nonlinearity $\tilde{\phi}$ is normalized to the componentwise sector $[-1,1]$, ensuring that the input and output satisfy the element-wise condition 
\begin{equation}
    \nu_{\phi,k}^{2}-x_{\phi,k}^{2}\geq0,
    \qquad
    k=1,\ldots,n_\phi,
    \quad
    \forall\nu_\phi\in[\underline\nu,\overline\nu]. \label{eq:normalized_sector} 
\end{equation}
\end{lemma}

\begin{pf}
Consider the NN controller \eqref{eq:nn_controller1}--\eqref{eq:nn_controller3} in feedback with the plant \eqref{eq:sys_trunc}, and let $(Z_{n_0}^*,u^*,\nu_\phi^*,\omega_\phi^*)$ be an equilibrium point. Define the deviation variables $\tilde Z_{n_0}(t) = Z_{n_0}(t)-Z_{n_0}^*$, $\tilde u(t)= u(t)-u^*$, $\tilde\nu_\phi(t) = \nu_\phi(t)-\nu_\phi^*$, $\tilde\omega_\phi(t) = \omega_\phi(t)-\omega_\phi^*$.

Subtracting the equilibrium relations from the affine representation \eqref{eq:nn_framework1} gives
\begin{align}
    \tilde u(t) &=N_{uz}\tilde Z_{n_0}(t)+N_{u\omega}\tilde\omega_\phi(t),\label{eq:equilibriumNN1}\\
    \tilde\nu_\phi(t) &=N_{\nu z}\tilde Z_{n_0}(t)+N_{\nu\omega}\tilde\omega_\phi(t).
    \label{eq:equilibriumNN2}
\end{align}

The corresponding translated activation function is
\begin{equation}
    \hat\phi(\tilde\nu_\phi(t)) \coloneqq \phi(\nu_\phi^*+\tilde\nu_\phi(t)) - \phi(\nu_\phi^*), \label{eq:translated_activation}
\end{equation}
so that $\tilde\omega_\phi(t)=\hat\phi(\tilde\nu_\phi(t))$, $\hat\phi(0)=0$.

For simplicity of notation, we henceforth drop the tildes and denote the translated nonlinearity $\hat\phi$ by $\phi$. The translated NN representation is therefore
\begin{align}
    \begin{bmatrix}
        u\\
        \nu_\phi
    \end{bmatrix} &= N\begin{bmatrix}
        Z_{n_0}\\
        \omega_\phi
    \end{bmatrix},\label{eq:nn_framework11}\\
    \omega_\phi(t) &= \phi(\nu_\phi(t)).\label{eq:nn_framework12}
\end{align}

The admissible set $\mathcal Z$ and the pre-activation interval $[\underline{\nu},\overline{\nu}]$ are henceforth understood in the translated coordinates.

For the fixed sector matrices $M_\phi$ and $R_\phi$, we apply the loop transformation shown in Figure~\ref{fig:NN-regularization} to the translated representation \eqref{eq:nn_framework11}--\eqref{eq:nn_framework12}. The transformation centers and scales the original sector so that the resulting nonlinearity belongs to the fixed componentwise sector $[-1,1]$. It consists of the following two steps:
\begin{enumerate}
    \item \textit{Centering}: Subtract the midpoint of the original sector, $c_\phi = \frac{M_{\phi}+R_{\phi}}{2}$, to define $\phi_{cent}(\nu_{\phi}(t)) = \phi(\nu_{\phi}(t)) - c_\phi \nu_{\phi}(t)$. By definition, $\phi_{cent}$ lies in the sector  $[-\rho_\phi, \rho_\phi]$, where $\rho_\phi = \frac{R_{\phi}-M_{\phi}}{2}$ represents the sector radius.
    \item \textit{Scaling}: Scale the centralized function by the inverse of the radius $\rho_\phi$, recovering \eqref{eq:tildephi_def}:
    \begin{equation} \label{eq:proof_tildephi}
        \hspace{-0.7cm}\tilde\phi(\nu_{\phi}(t)) = \rho_\phi^{-1}\phi_{cent}(\nu_{\phi}(t)) = \rho_\phi^{-1}(\phi(\nu_{\phi}(t)) - c_\phi\nu_{\phi}(t)).
    \end{equation}
\end{enumerate}

It follows that
\begin{equation*}
    |\tilde\phi_i(\nu_{\phi,i}(t))|\leq|\nu_{\phi,i}(t)|,\qquad i\in\{1,\dots,n_\phi\},
\end{equation*}
which establishes the normalized componentwise sector condition. \redx{Equivalently, \(\nu_{\phi,k}^2-x_{\phi,k}^2\geq0\) for every \(k\), which proves the componentwise condition \eqref{eq:normalized_sector}.}

Rearranging \eqref{eq:proof_tildephi} gives the loop transformation \eqref{eq:loop_transform}:
\begin{equation} \label{eq:proof_wphi}
    \omega_{\phi}(t) = \rho_\phi x_{\phi}(t)+c_\phi \nu_{\phi}(t),
\end{equation}
\noindent where $x_{\phi}(t) = \tilde{\phi}(\nu_{\phi}(t))$ and $\omega_{\phi}(t)=\phi(\nu_{\phi}(t))$. Substituting \eqref{eq:proof_wphi} into the affine relation \eqref{eq:nn_framework11}, and using the definitions of $C_1,\dots,C_4$ given in the statement of the lemma, yields:
\begin{align}
    u(t) &= N_{uz} Z_{n_{0}}(t) + C_1 x_{\phi}(t) + C_2 \nu_{\phi}(t), \label{eq:u_n} \\
    \nu_{\phi}(t) &= N_{\nu z} Z_{n_{0}}(t) + C_3 x_{\phi}(t) + C_4 \nu_{\phi}(t). \label{eq:v_phi}
\end{align}

Since $(I_{n_\phi} - C_4)$ is invertible, \eqref{eq:v_phi} can be solved explicitly for the activation input $\nu_{\phi}(t)$:
\begin{align}
    \hspace{-0.7cm}\nu_{\phi}(t) = (I_{n_\phi} - C_4)^{-1} N_{\nu z} Z_{n_{0}}(t) + (I_{n_\phi} - C_4)^{-1} C_3 x_{\phi}(t). \label{eq:v_phi2}
\end{align}

Substituting \eqref{eq:v_phi2} into \eqref{eq:u_n} produces the closed-form expression for $u(t)$:
\begin{align}
    u(t) &= \left(N_{uz} + C_2 (I_{n_\phi} - C_4)^{-1} N_{\nu z}\right) Z_{n_{0}}(t) \nonumber\\
    &\quad+\left(C_1 + C_2 (I_{n_\phi} - C_4)^{-1} C_3\right) x_{\phi}(t). \label{eq:u_n2}
\end{align}

Concatenating \eqref{eq:v_phi2} and \eqref{eq:u_n2} recovers exactly the matrix $\tilde{N}$ defined in \eqref{eq:nn_normalized}, establishing the equivalent normalized linear mapping $\mathcal{F}: N \mapsto \tilde{N}$. This formulation isolates the nonlinearity $\tilde\phi$ within the normalized componentwise sector $[-1,1]$, subject to the componentwise QCs \eqref{eq:normalized_sector}.

This concludes the proof. \hfill $\square$
\end{pf}

\section{Stability Analysis of the Infinite- Dimensional Closed-Loop System}\label{section:stability_analysis}

The boundary-actuation coefficients satisfy ~$|\beta_n| \allowbreak=\allowbreak O(\sqrt{\lambda_n})$ and are therefore not guaranteed to be square summable. Nevertheless, Lemma~\ref{lem:effective_spillover} establishes that $\left\{\beta_n/\lambda_n\right\} _{n\geq1}\allowbreak\in\ell^2$.

This high-frequency behavior motivates a composite Lyapunov functional with two terms. The first term measures the retained state $Z_{n_0}$ used by the NN controller, whereas the second measures the residual modal coordinates with the spectral weights $1/\lambda_n$. These weights are consistent with the negative-order state space naturally associated with the unbounded boundary input and permit the residual input contribution to be controlled through the summable quantity $\displaystyle S_2(n_0) = \textstyle\sum_{n=n_0+1}^{\infty} \beta_n^2/\lambda_n^2$. This construction is used below to certify the complete infinite-dimensional closed-loop system \citep{Lhachemi26}.

To account for this high-frequency behavior, we introduce a Lyapunov functional comprising two terms. The first term accounts for the first $n_0$ modes through the truncated state $Z_{n_0}$ in \eqref{eq:sys_trunc}, which is also used for the expert-controller design. The second term accounts for the residual modes in the modal coordinates $z_n$ of \eqref{eq:modes}, weighted by $1/\lambda_n$. This weighting yields the square-summable tail coefficient sequence $\left\{\beta_n/\lambda_n\right\}_{n>n_0}\in\ell^2$, since $\sum_{n>n_0}\beta_n^2/\lambda_n^2<\infty$, as established in Lemma~\ref{lem:effective_spillover}. The resulting Lyapunov functional is coercive with respect to a norm equivalent to the norm of $\mathcal H_{-1/2}$ and provides a stability certificate for the full infinite-dimensional closed loop \citep{Lhachemi26}.

\begin{remark}\label{rem:star_norm}
Let $\mathcal A$ be the positive self-adjoint Sturm\-/Liouville operator with eigenpairs $(\lambda_n,\Phi_n)_{n\geq1}$. Define $\mathcal H_{1/2}\coloneqq D(\mathcal A^{1/2})$, $\mathcal H_{-1/2}\coloneqq \mathcal H_{1/2}'$, where the duality is taken with respect to the pivot space $L^2(0,1)$ \cite[Section~3.4]{tucsnak2009observation}. The space $\mathcal H_{-1/2}$ is endowed with the norm
\begin{equation}
    \|z\|_{\mathcal H_{-1/2}}^2\coloneqq\|\mathcal A^{-1/2}z\|_{L^2(0,1)}^2 =\sum_{n=1}^{\infty}\frac{|z_n|^2}{\lambda_n}.
    \label{eq:Hminus_norm}
\end{equation}

For a fixed truncation order $n_0$, define
\begin{equation}
    \|z\|_\star^2\coloneqq\sum_{n=1}^{n_0}|z_n|^2+\sum_{n=n_0+1}^{\infty}\frac{|z_n|^2}{\lambda_n}.\label{eq:star_norm}
\end{equation}

Equivalently,
\[
    \|z\|_\star^2 =\|Z_{n_0}\|^2+\sum_{n=n_0+1}^{\infty}\frac{|z_n|^2}{\lambda_n}.
\]

Since $\lambda_1\leq\lambda_n\leq\lambda_{n_0}$ for $n=1,\ldots,n_0$, it follows that
\begin{equation*}
   \hspace{-0.5cm} \min\{1,\lambda_1\}\|z\|_{\mathcal H_{-1/2}}^2\leq \|z\|_\star^2\leq\max\{1,\lambda_{n_0}\}\|z\|_{\mathcal H_{-1/2}}^2.
\end{equation*}

Hence, $\|\cdot\|_\star$ is an equivalent norm on $\mathcal H_{-1/2}$, with equivalence constants depending on the fixed truncation order $n_0$. This norm reflects the structure of the Lyapunov functional: the retained modes are measured in the Euclidean norm, whereas the residual modes are measured in the negative-order norm generated by $\mathcal A$.
\end{remark}

\begin{theorem} \label{theo}
Consider the reaction--diffusion system described by \eqref{eq:pde_main}--\eqref{eq:pde_ini}, with $p \in C^{1}([0,1])$ such that $p(x)>0$, $\forall x\in[0,\,1]$, and $q_{c} \in \mathbb{R}$. Let $\{(\lambda_n,\Phi_n)\}_{n\ge1}$ denote the associated Sturm--Liouville eigenpairs, and let $\{z_n\}_{n\ge1}$ be the corresponding modal coordinates satisfying \eqref{eq:modes}. Let $n_{0} \geq 1$ and $\delta > 0$ be such that
\begin{equation*}
    -\lambda_{n} + q_{c} < -\delta < 0, \qquad \forall n \geq n_{0}+1.
\end{equation*}

Let $\mathcal Z_h \coloneqq\left\{Z\in\mathbb R^{n_0}:\underline Z\leq Z\leq\overline Z\right\}$ be the hyperrectangle used to compute the sector bounds, and let $m_\phi,r_\phi\in\mathbb R^{n_\phi}$ be the corresponding componentwise sector-bound vectors computed according to Lemma~\ref{lemma:sector_condition}. Consider the NN controller \eqref{eq:nn_framework1}--\eqref{eq:nn_framework2} with equilibrium point $Z_{n_{0}}^{*} = 0$ and state constraint set $\mathcal{Z}$ defined in \eqref{eq:admissible_set}, and let $\tilde N$ denote its loop-transformed representation defined in Lemma~\ref{lemma:loop_transformation}.

Assume that $\mathcal Z\subseteq\mathcal Z_h$. Fix scalars $\gamma>0$ and $\mu\geq1/2$. Suppose that there exist a positive definite matrix $P\in\mathbb S_{++}^{n_0}$ and a vector $\theta\in\mathbb R^{n_\phi}$ satisfying $\theta_i>0$ for every $i\in\{1,\dots,n_\phi\}$, such that the following inequalities hold:
\begin{align}
    \hspace{-0.9cm}\tilde{R}_{\nu}^{\top} \bar{Q}_{1} \tilde{R}_{\nu}
    + \tilde{R}_{\phi}^{\top}Q_{2}\tilde{R}_{\phi}
    + \mathcal{R} &\prec 0, \label{eq:lmi} \\
    \Gamma_{n_{0}+1} & \leq 0, \label{eq:lmi2} \\
    \begin{bmatrix} s_{i}^{2} & S_{i}^{\top} \\ S_{i} & P \end{bmatrix}
    &\succeq 0, \quad i=1,\dots,n_{z}, \label{eq:lmi4}
\end{align}
\noindent where $S_{i}^{\top}$ denotes the $i$-th row of $S$, and the first-tail coefficient $\Gamma_{n_0+1}$ and the spillover matrix $\mathcal R$ are defined by
\begin{align}
    \Gamma_{n_0+1} &\coloneqq \left(2\left(\frac{q_{c}+\delta}{\lambda_{n_{0}+1}} - 1\right) + \frac{1}{\mu}\right), \label{eq:Gamma} \\
    \mathcal{R} &\coloneqq \gamma\mu S_{2}(n_{0}) \tilde{R}_{u}^{\top}\tilde{R}_{u}, \qquad S_{2}(n_{0})\coloneqq \sum_{n\geq n_{0}+1}\frac{\beta_{n}^{2}}{\lambda_{n}^{2}}. \label{eq:spill}
\end{align}

The remaining finite-dimensional block matrices are given by $\Theta = \diag(\theta)$ and
\begin{align*}
    \bar{Q}_{1} &\coloneqq \begin{bmatrix}
    A_{n_{0}}^{\top}P + PA_{n_{0}} + 2\delta P & PB_{n_{0}} \\
    B_{n_{0}}^{\top}P & 0
    \end{bmatrix},\\
        Q_{2} &\coloneqq \begin{bmatrix}
        \Theta & 0 \\
        0 & -\Theta
    \end{bmatrix},
    \\
    \tilde{R}_{\nu} &\coloneqq \begin{bmatrix}
        I_{n_{0}} & 0 \\
        \tilde{N}_{u z} & \tilde{N}_{u x}
    \end{bmatrix},
    \\
    \tilde{R}_{\phi} &\coloneqq \begin{bmatrix}
        \tilde{N}_{\nu z} & \tilde{N}_{\nu x} \\
        0 & I_{n_{\phi}}
    \end{bmatrix},
    \\
    \tilde{R}_{u} &\coloneqq \begin{bmatrix} \tilde{N}_{uz} & \tilde{N}_{ux} \end{bmatrix}.
\end{align*}

Therefore, the origin of the closed-loop system \eqref{eq:modes} is locally exponentially stable in $\mathcal H_{-1/2}$, or equivalently with respect to the norm $\|\cdot\|_\star$ defined in \eqref{eq:star_norm}. Moreover, the set 
\begin{align}
    \hspace{-0.85cm} \mathcal{E}(P,\gamma) \coloneqq \Bigg\{ z \in \mathcal H_{-1/2} :
     Z_{n_{0}}^{\top}PZ_{n_{0}} + \gamma\!\! \sum_{n\geq n_{0}+1} \frac{z_{n}^{2}}{\lambda_{n}} \leq 1 \Bigg\}, \label{eq:set_invariance}
\end{align}
is positively invariant and constitutes an inner approximation of the region of attraction in $\mathcal H_{-1/2}$. In particular, there exists a constant $M\geq1$, independent of the initial condition, such that, for every $z(0)\in\mathcal E(P,\gamma)$, the corresponding mild solution globally exists and satisfies
\begin{equation*}
    \|z(t,\cdot)\|_\star \le M e^{-\delta t} \|z(0,\cdot)\|_\star, \qquad \forall t\ge 0.
\end{equation*}
\end{theorem}

\begin{pf}
By Lemma~\ref{lemma:lipschitz_nn}, the tanh NN policy is globally Lipschitz and globally bounded. Together with Lemma~\ref{lem:effective_spillover}, this verifies the hypotheses of Lemma~\ref{lem:closed_loop_wellposedness}. Hence the closed-loop system admits a unique global mild solution $z\in C\bigl([0,\infty);\mathcal H_{-1/2}\bigr)$. The infinite-dimensional Lyapunov calculations below are understood through finite modal sums followed by passage to the limit.


Consider the Lyapunov functional
\begin{equation}
    V(t) = V_{1}(t) + V_{2}(t),\label{eq:lyap_functional}
\end{equation}
\noindent where\footnote{By Lemma~\ref{lem:closed_loop_wellposedness}, $z(t,\cdot)\in L^2(0,1)$ for every $t>0$. More rigorously, the calculations are first performed on $[\varepsilon,T]$ for the partial sums $V_2^{(M)}(t)=\gamma\sum_{n=n_0+1}^{M}z_n(t)^2/\lambda_n$. The limit $M\to\infty$ is then taken using the corresponding summability estimates, followed by $\varepsilon\to0$ using continuity in $\mathcal H_{-1/2}$.}
\begin{align}
    V_{1}(t) &= Z_{n_{0}}^{\top}(t)P Z_{n_{0}}(t), \label{eq:V1} \\
    V_{2}(t) &= \gamma \sum_{n\geq n_{0}+1}\frac{z_{n}^{2}(t)}{\lambda_{n}}\label{eq:V2} ,
\end{align}
\noindent and $P\in\mathbb S_{++}^{n_0}$ is positive definite and $\gamma>0$ is a fixed scalar parameter.

By the definition of the norm $\|\cdot\|_\star$ introduced in
Remark~\ref{rem:star_norm},
\begin{equation*}
    \|z(t,\cdot)\|_\star^2 = \|Z_{n_0}(t)\|^2 + \sum_{n>n_0}\frac{z_n(t)^2}{\lambda_n}.
\end{equation*}

Since $P\succ 0$, $\underline c\,\|z(t,\cdot)\|_\star^2 \le V(t) \le \overline c\,\|z(t,\cdot)\|_\star^2$, where $\underline c = \min\{\lambda_{\min}(P),\gamma\}$, $\overline c = \max\{\lambda_{\max}(P),\gamma\}$. Hence, the Lyapunov functional is coercive and continuous with respect to $\|\cdot\|_\star$. By Remark~\ref{rem:star_norm}, $\|\cdot\|_\star$ is an equivalent norm on $\mathcal H_{-1/2}$. Consequently, $V$ is also coercive and continuous with respect to the canonical norm of $\mathcal H_{-1/2}$.

Evaluating the time derivative of \eqref{eq:V1} along the trajectories of the truncated system \eqref{eq:sys_trunc}, yields
\begin{align}
   \hspace{-0.9cm} \dot{V}_{1}(t) &= 
     \begin{bmatrix}
        Z_{n_{0}}(t) \\ u(t)
    \end{bmatrix}^{\!\top} \begin{bmatrix}
        A_{n_{0}}^{\top}P + PA_{n_{0}} & PB_{n_{0}} \\ B_{n_{0}}^{\top}P & 0 
    \end{bmatrix} \begin{bmatrix}
        Z_{n_{0}}(t) \\ u(t)
    \end{bmatrix}. \label{eq:dV1}
\end{align}

By assumption, the closed-loop equilibrium is $Z_{n_0}^\ast=0$. Moreover, the NN controller satisfies $\pi(0)=0$, and therefore the corresponding equilibrium input is $u^\ast=0$. Based on Lemmas~\ref{lemma:sector_condition} and \ref{lemma:loop_transformation}, define
\begin{align*}
    Q_{1} \coloneqq \begin{bmatrix}
    A_{n_{0}}^{\top}P + PA_{n_{0}} & PB_{n_{0}} \\ B_{n_{0}}^{\top}P & 0 
    \end{bmatrix},
\end{align*}
\noindent and then plugging $Q_{1}$ and \eqref{eq:u_n2} in \eqref{eq:dV1}, it follows that $\dot{V}_{1}$ can be rewritten as
\begin{align}
    \dot{V}_{1}(t) &= 
     \begin{bmatrix} Z_{n_{0}}(t) \\ x_{\phi}(t) \end{bmatrix}^{\!\top} \tilde{R}_{\nu}^{\top} Q_{1}\tilde{R}_{\nu} \begin{bmatrix}
        Z_{n_{0}}(t) \\ x_{\phi}(t)
    \end{bmatrix}. \label{eq:dotV1}
\end{align}

To incorporate the prescribed exponential decay rate $\delta$ into the finite-dimensional Lyapunov estimate, define the aggregated state $\xi(t)\coloneqq \begin{bmatrix} Z_{n_0}(t)^\top & x_\phi(t)^\top \end{bmatrix}^\top$. From \eqref{eq:dotV1} and the definition of $\bar Q_1$, it follows
that
\begin{equation}
    \dot V_1(t)+2\delta V_1(t)=\xi(t)^\top\tilde R_\nu^\top\bar Q_1\tilde R_\nu\xi(t).\label{eq:dotV1_decay}
\end{equation}

The condition
\begin{align*}
    \xi(t)^{\top} \tilde R_\nu^{\top} \bar{Q}_{1} \tilde R_\nu \, \xi(t)\leq 0,
\end{align*}
is equivalent to $\dot V_1(t)+2\delta V_1(t)\le0$, since $\bar{Q}_{1}$ augments $Q_{1}$ with the term $2\delta P$, as follows from \eqref{eq:dotV1}.


For every time $t$ such that $Z_{n_0}(t)\in\mathcal Z$, the inclusion $\mathcal Z\subseteq\mathcal Z_h$ and Lemma~\ref{lemma:sector_condition} guarantee that each activation input-output pair satisfies its corresponding local sector condition. By the loop transformation of Lemma~\ref{lemma:loop_transformation}, these inequalities become \eqref{eq:normalized_sector}. Using the diagonal multiplier matrix $\Theta=\operatorname{diag}(\theta_1,\ldots,\theta_{n_\phi})\succ0$, multiplying the $k$-th inequality by $\theta_k>0$ and summing over $k=1,\ldots,n_\phi$ yields
\begin{equation*}
    \begin{bmatrix}
        \nu_\phi(t)\\ x_\phi(t)
    \end{bmatrix}^{\!\top}
    \begin{bmatrix}
        \Theta & 0\\
        0 & -\Theta
    \end{bmatrix}
    \begin{bmatrix}
        \nu_\phi(t)\\ x_\phi(t)
    \end{bmatrix}
    \geq0.
\end{equation*}

Substituting \eqref{eq:v_phi2} into this weighted quadratic constraint yields
\begin{align}
    \xi(t)^{\!\top} \tilde{R}_{\phi}^{\top}Q_{2}\tilde{R}_{\phi} \,\xi(t) &\geq 0.\label{eq:qc_sector_cond_final}
\end{align}

Since the weighted quadratic constraint \eqref{eq:qc_sector_cond_final} is nonnegative, adding it to \eqref{eq:dotV1_decay} yields
\begin{align}
    \hspace{-0.7cm}\dot{V}_{1}(t) + 2\delta V_{1}(t) \le \xi(t)^{\!\top} \Big( \tilde{R}_{\nu}^{\top} \bar{Q}_{1} \tilde{R}_{\nu} + \tilde{R}_{\phi}^{\top} Q_{2} \tilde R_\phi \Big) \xi(t). \label{eq:dotV1_bound}
\end{align}

We now consider the derivative of $V_2$. Differentiating \eqref{eq:V2} with respect to time and substituting \eqref{eq:modes}, we obtain
\begin{align}
   \hspace{-0.8cm} \dot{V}_{2}(t) = 2\gamma \sum_{n\geq n_{0}+1} z_{n}(t)\left( \left(\frac{q_{c}}{\lambda_{n}} - 1\right)z_{n}(t) + \frac{\beta_{n}}{\lambda_{n}}u (t) \right). \label{eq:dotV2}
\end{align}

Applying Young's inequality parameterized by $\mu \geq 1/2$ to decouple the cross-terms in \eqref{eq:dotV2}, we obtain the following inequality:
\begin{align}
    \hspace{-0.75cm} 2 \sum_{n\geq n_{0}+1} \frac{\beta_{n}}{\lambda_{n}} z_{n}(t) u(t)
    &\leq \frac{1}{\mu} \sum_{n\geq n_{0}+1} z_{n}^{2}(t)
    + \mu S_{2}(n_{0}) u^{2}(t),\label{eq:residual_norms}
\end{align}
adding the  decay margin $2\delta V_{2}(t)$ yields the following inequality:
\begin{align*}
    \hspace{-0.85cm}\dot{V}_{2}(t) + 2\delta V_{2}(t) &\leq \gamma \sum_{n\geq n_{0}+1} \left(2  \left( \frac{q_{c}+\delta}{\lambda_{n}} -1 \right) + \frac{1}{\mu} \right) z_{n}^2(t) \\
    &\quad +\gamma\mu S_{2}(n_{0})u^2(t).
\end{align*}

From the normalized controller representation \eqref{eq:nn_tilde_framework1} in Lemma~\ref{lemma:loop_transformation} and the definition of $\tilde{R}_{u}$, we have  $$u(t)^2 = \xi(t)^{\top} \tilde R_u^{\top} \tilde R_u \,\xi(t).$$ 

Then, it follows that
\begin{align}
    \hspace{-0.85cm}\dot{V}_{2}(t) + 2\delta V_{2}(t) &\leq \gamma \sum_{n\geq n_{0}+1} \left(2 \left( \frac{q_{c}+\delta}{\lambda_{n}} - 1 \right)+ \frac{1}{\mu} \right) z_{n}^2(t) \nonumber\\
    &\quad+\xi(t)^{\top}\mathcal{R}\, \xi(t). \label{eq:dotV2_bound}
\end{align}

Combining \eqref{eq:dotV1_bound} and \eqref{eq:dotV2_bound} yields the following bound for the derivative of the
full Lyapunov functional:
\begin{align}
    \dot V(t) + 2\delta V(t) &\le \xi(t)^{\!\top} \Big(\tilde R_\nu^\top \bar{Q}_{1} \tilde R_\nu + \tilde R_\phi^\top Q_2 \tilde R_\phi + \mathcal{R} \Big) \xi(t) \nonumber \\
    &\quad + \gamma \sum_{n\geq n_{0}+1} \left(2  \left( \frac{q_{c}+\delta}{\lambda_{n}} - 1 \right)+ \frac{1}{\mu} \right) z_{n}^2(t). \label{eq:V_full1}
\end{align}

For every $n\geq n_0+1$, define
\begin{equation}
    \Gamma_n\coloneqq 2\left(\frac{q_c+\delta}{\lambda_n}-1\right)+ \frac{1}{\mu}.
    \label{eq:Gamma_n}
\end{equation}

We now show that $\Gamma_n\leq0$ for every $n\geq n_0+1$. Since $\{\lambda_n\}_{n\geq1}$ is strictly increasing, the monotonicity of $\Gamma_n$ depends on the sign of $q_c+\delta$.

If $q_c+\delta\geq0$, then $\Gamma_n$ is nonincreasing, and therefore $\Gamma_n\leq \Gamma_{n_0+1}\leq0$, $n\geq n_0+1$, where the last inequality follows from \eqref{eq:lmi2}. If $q_c+\delta<0$, then $\frac{q_c+\delta}{\lambda_n}<0$, and hence $\Gamma_n < -2+\frac{1}{\mu}\leq 0$, $n\geq n_0+1$, where the last inequality follows from $\mu\geq1/2$.

Consequently, $\Gamma_n\leq0$, $\forall n\geq n_0+1$, and thus $$\gamma \sum_{n\geq n_0+1} \Gamma_nz_n^2(t) \leq0.$$

It follows from \eqref{eq:V_full1} that
\begin{align*}
    \hspace{-0.7cm}\dot V(t)+2\delta V(t) &\leq\xi(t)^\top \Big(\tilde R_\nu^\top\bar Q_1\tilde R_\nu + \tilde R_\phi^\top Q_2\tilde R_\phi+\mathcal R\Big)\xi(t).
\end{align*}

Finally, condition \eqref{eq:lmi} implies
\begin{equation}
    \dot V(t)+2\delta V(t) \leq 0.
    \label{eq:V_decay}
\end{equation}

We now prove the positive invariance of $\mathcal{E}(P,\gamma)$ and the regional exponential stability claim. Since the quadratic sector constraints bounding the NN are local, define
\begin{equation*}
    T^\star \coloneqq \sup\left\{ T>0:\; V(t)\leq1,\quad\forall t\in[0,T) \right\}.
\end{equation*}

For every $t\in[0,T^\star)$, we have $Z_{n_0}(t)^\top PZ_{n_0}(t) \leq V(t)\leq1$. Thus,
\[
    Z_{n_0}(t)\in\mathcal E_P,
    \qquad
    \mathcal E_P
    \coloneqq
    \left\{
        Z_{n_0}\in\mathbb R^{n_0}:
        Z_{n_0}^\top PZ_{n_0}\leq1
    \right\}.
\]

By the Schur complement, condition \eqref{eq:lmi4} implies
\begin{equation*}
    S_i^\top P^{-1}S_i\leq s_i^2,\qquad i\in\{1,\dots,n_z\}.
\end{equation*}

Hence, for every $Z_{n_0}\in\mathcal E_P$,
\begin{align*}
    |S_i^\top Z_{n_0}|^2\leq\bigl(Z_{n_0}^\top PZ_{n_0}\bigr)\bigl(S_i^\top P^{-1}S_i\bigr)\leq s_i^2.
\end{align*}

Therefore, $\mathcal E_P\subseteq\mathcal Z$. 

Since $\mathcal Z\subseteq\mathcal Z_h$, the local sector conditions of Lemma~\ref{lemma:sector_condition} are valid on $[0,T^\star)$. Therefore, the feasibility of \eqref{eq:lmi} and \eqref{eq:lmi2} implies
\[
    \dot V(t)+2\delta V(t)\leq0,
    \qquad
    \forall t\in[0,T^\star).
\]

Integrating this differential inequality gives
\begin{equation}
    V(t) \leq e^{-2\delta t}V(0),\qquad\forall t\in[0,T^\star).\label{eq:lyap_bound}
\end{equation}

Let $z(0)\in\mathcal E(P,\gamma)$, so that $V(0)\leq1$. It follows from \eqref{eq:lyap_bound} that
\[
    V(t)\leq e^{-2\delta t}V(0)\leq1,\qquad \forall t\in[0,T^\star).
\]

The estimate above and the continuity of $V$, guaranteed by Lemma~\ref{lem:closed_loop_wellposedness}, exclude any finite exit from the sublevel set $\{z:V(z)\leq1\}$. Hence, $T^\star=\infty$. Consequently,
\[
    V(t)\leq e^{-2\delta t}V(0), \qquad \forall t\geq0,
\]
and $\mathcal E(P,\gamma)$ is positively invariant.

Finally, the norm equivalence gives
\begin{align*}
    \|z(t,\cdot)\|_\star^2 \leq \frac{1}{\underline c}V(t)\leq \frac{\overline c}{\underline c}e^{-2\delta t}
    \|z(0,\cdot)\|_\star^2.
\end{align*}
Taking square roots yields 
\begin{equation*}
    \|z(t,\cdot)\|_\star\leq M e^{-\delta t}\|z(0,\cdot)\|_\star,
\end{equation*}
where $M\coloneqq\sqrt{\frac{\overline c}{\underline c}}\geq1$.

This proves that every solution initiating in $\mathcal{E}(P,\gamma)$ exists globally, remains within $\mathcal{E}(P,\gamma)$ for all $t \ge 0$, and exponentially converges to the origin with respect to the norm $\|\cdot\|_\star$, or equivalently in $\mathcal H_{-1/2}$.\hfill $\square$
\end{pf}

\begin{remark}
The space $\mathcal H_{-1/2}$ is weaker than $L^2(0,1)$ because its norm attenuates the high-frequency modal coordinates by $1/\lambda_n$. Therefore, Theorem~\ref{theo} does not claim uniform $L^2$ exponential stability for arbitrary initial conditions in $\mathcal H_{-1/2}$.
\end{remark}

{\color{black}

The feasibility conditions of Theorem~\ref{theo} are not convex in the NN weights and the Lyapunov variables jointly. Indeed, the loop-transformed controller matrix $\tilde N=\mathcal F(N)$ depends nonlinearly on the NN weights $N$ and appears bilinearly coupled with the Lyapunov variables.

The scalar parameters $\gamma>0$ and $\mu\ge1/2$ are fixed before solving the training problem. In particular, $\mu$ is selected so that the residual-mode condition $\Gamma_{n_0+1}\le0$ is satisfied. Therefore, neither $\gamma$ nor $\mu$ is treated as a decision variable in the semidefinite training subproblem.

For a fixed NN controller $N$, the matrix $\tilde N=\mathcal F(N)$ is known. Consequently, $\tilde R_\nu$, $\tilde R_\phi$, and $\tilde R_u$ are known, and the direct-spillover stability condition \eqref{eq:lmi} is convex in the Lyapunov variables $P$ and $\Theta$.

To derive a convex representation suitable for the certificate\-/update subproblem, we first consider the spillover-free part of \eqref{eq:lmi}. Although the inequality is still nonlinear in the variables $P$, $\Theta$ and $\tilde{N}$, applying the Schur complement gives
\begin{align*}
    \begin{bmatrix} -\Xi_P & -P B_{n_{0}} \tilde{N}_{ux} & \tilde{N}_{\nu z}^{\top} \\ -\tilde{N}_{ux}^{\top} B_{n_{0}}^{\top} P & \Theta & \tilde{N}_{\nu x}^{\top} \\ \tilde{N}_{\nu z} & \tilde{N}_{\nu x} & \Theta^{-1} \end{bmatrix} \succ 0,
\end{align*}
\noindent where $\Xi_P \coloneqq \big(A_{n_{0}}^{\top} P + P A_{n_{0}} + P B_{n_{0}} \tilde{N}_{uz} + \tilde{N}_{uz}^{\top} B_{n_{0}}^{\top} P \,\allowbreak + 2\delta P \big)$.

To solve the bilinear coupling, we pre- and post-multiply the previous inequality by $\diag(P^{-1}, \Theta^{-1}, I_{n_{\phi}})$. Defining 
\begin{align}
    H_{1} &\coloneqq P^{-1} \in \mathbb{S}_{++}^{n_0},
    &
    H_{2} &\coloneqq \Theta^{-1} \in \mathbb{D}_{++}^{n_\phi}, \label{def:1}
    \\
    L_{1} &\coloneqq \tilde{N}_{uz}H_{1}
    \in \mathbb{R}^{n_u \times n_0},
    &
    L_{2} &\coloneqq \tilde{N}_{ux}H_{2}
    \in \mathbb{R}^{n_u\times n_\phi}, \label{def:2}
    \\
    L_{3} &\coloneqq \tilde{N}_{\nu z}H_{1}
    \in \mathbb{R}^{n_\phi\times n_0},
    &
    L_{4} &\coloneqq \tilde{N}_{\nu x}H_{2}
    \in \mathbb{R}^{n_\phi\times n_\phi}. \label{def:3}
\end{align}

Moreover, let
\begin{equation*}
    M_H(\bar H,\bar L) \coloneqq
    \begin{bmatrix}
        -\Xi_H &-B_{n_0}L_2 &L_3^\top\\
        -L_2^\top B_{n_0}^\top &H_2 &L_4^\top\\
        L_3 &L_4 &H_2
    \end{bmatrix}, 
\end{equation*}
where
\begin{equation*}
    \Xi_H \coloneqq H_1A_{n_0}^\top+A_{n_0}H_1 + B_{n_0}L_1+L_1^\top B_{n_0}^\top + 2\delta H_1.
\end{equation*}

The spillover term can then be included directly by defining $L_u \coloneqq \begin{bmatrix} L_1 & L_2 & 0 \end{bmatrix} \in \mathbb{R}^{n_u\times(n_0+2n_\phi)}$, and hence applying a second Schur complement:
\begin{equation*}
    \mathcal M_H(\bar H,\bar L)
    \coloneqq
    \begin{bmatrix}
        M_H(\bar H,\bar L)
        &
        \sqrt{\gamma\mu S_2(n_0)}\,L_u^\top
        \\[1mm]
        \sqrt{\gamma\mu S_2(n_0)}\,L_u
        &
        I_{n_u}
    \end{bmatrix},
\end{equation*}
where for a scalar boundary input, $n_u=1$ and the lower-right block reduces to the scalar $1$.

Indeed, by the Schur complement $\mathcal M_H(\bar H,\bar L)\succ0$ is equivalent to
\begin{equation}
    \mathcal M_H(\bar H,\bar L) \succ0
    \quad\Longleftrightarrow\quad
    M_H(\bar H,\bar L)
    -\gamma\mu S_2(n_0)L_u^\top L_u
    \succ0.
    \label{eq:schur_spillover}
\end{equation}

Following the aggregate notation used in \cite{yin2021imitation}, define
\begin{equation}
    \operatorname{LMI}(\bar H,\bar L)\succ0
    \quad\Longleftrightarrow\quad
    \left\{
    \begin{aligned}
    &\mathcal M_H(\bar H,\bar L)\succ0,\\[1mm]
    &
    \begin{bmatrix}
    s_i^2&S_i^\top H_1\\
    H_1S_i&H_1
    \end{bmatrix}
    \succeq0,
    \\[1mm]
    &H_1\succ0,\qquad H_2\succ0,
    \end{aligned}
    \right.
    \label{eq:training_LMI}
    \end{equation}
where $i\in\{1,\,\ldots,\,n_z\}$, $\bar H=\{H_1,\,H_2\}$, $\bar L=\{L_1,\,L_2,\,L_3,\,\allowbreak L_4\}$.




The residual terms arising from the infinite-dimensional tail are directly included in the training problem. Although the tail coefficient $S_2(n_0)$ depends only on the PDE and the truncation order, the spillover matrix $\mathcal R = \gamma\mu S_2(n_0) \tilde R_u^\top\tilde R_u$ also depends on the controller variables. For fixed values of $\gamma$ and $\mu$, the second Schur complement above provides an equivalent representation that is affine in the transformed decision variables $\bar H$ and $\bar L$.

Therefore, the aggregate LMI \eqref{eq:training_LMI} incorporates the finite-dimensional stability, state-set containment, and direct-spillover requirements, while the scalar residual-tail condition \eqref{eq:lmi2} is enforced separately through the fixed choice of $\mu$. After the algorithm convergence, the same condition is recomputed using the fixed NN controller, its exact loop-transformed representation, and the corresponding local sector bounds. This final computation validates the obtained controller but does not constitute a separate spillover optimization problem.

For a prescribed exponential decay rate $\delta>0$ and a selected truncation order $n_0$, the parameters $\gamma>0$ and $\mu\geq1/2$ are fixed before training, with $\mu$ chosen so that $\Gamma_{n_0+1}\leq0$. The sufficient conditions of Theorem~\ref{theo} can then be imposed through the single semidefinite constraint \eqref{eq:training_LMI}. Increasing $n_0$ reduces the residual coefficient $S_2(n_0)$ and may therefore facilitate feasibility of the direct-spillover training LMI. Nevertheless, feasibility must be checked for each selected truncation order and is not guaranteed solely by increasing $n_0$.
}

\subsection{Feasibility of the Constrained Learning Problem}

The constrained learning framework incorporates finite-dimensional stability, and infinite-dimensional spillover robustness into the single aggregate training condition \eqref{eq:training_LMI}. The residual-mode requirement is enforced by fixing $\mu\geq1/2$ so that $\Gamma_{n_0+1}\leq0$ before solving the training problem.

To establish feasibility for all sufficiently large truncation orders, we require that the spillover-free finite-dimensional approximation admit a \emph{uniformly robust Lyapunov certificate}. Its non-zero stability margin provides the mechanism needed to absorb the direct-spillover term in \eqref{eq:training_LMI} as $S_2(n_0)$ tends to zero.

To study feasibility of the direct-spillover training condition as the truncation order increases, we formulate the required uniform robustness property directly in the analysis variables $(P,\Theta)$. This property concerns the nominal finite-dimensional dissipation and is used only to establish that the full spillover-dependent condition becomes feasible for all sufficiently large $n_0$. A uniform bound on the transformed gain $\tilde R_u=\begin{bmatrix}\tilde N_{uz} & \tilde N_{ux}\end{bmatrix}$ is likewise required, since the controller is retrained for each truncation order.

\begin{assumption} \label{ass:certifiable_robust_expert}
For a prescribed decay rate $\delta>0$, there exist constants $\eta>0$, $\kappa_u>0$ and $N_0\in\mathbb{N}$ such that, for every truncation order $n_0 \ge N_0$, the finite-dimensional closed-loop system induced by the candidate NN controller in the chosen architecture admits a robust Lyapunov certificate that is uniform in $n_0$. Specifically, for each $n_0 \ge N_0$ there exist $P \succ 0$, $\Theta = \diag(\theta) \succ 0$, and a transformed representation $\tilde N$ associated with a candidate NN controller in the chosen architecture such that the controller is required to satisfy the componentwise local sector conditions of Lemma~\ref{lemma:sector_condition} and the containment condition \eqref{eq:lmi4}, and
\begin{equation}
    \tilde{R}_{\nu}^{\top} \bar{Q}_{1} \tilde{R}_{\nu}
    + \tilde{R}_{\phi}^{\top} Q_{2} \tilde{R}_{\phi}
    \preceq -\eta I,
    \qquad
    \|\tilde R_u\| \le \kappa_u,
    \label{eq:uniform_margin}
\end{equation}
\noindent where the margin $\eta$ and the gain bound $\kappa_u$ are independent of the truncation order $n_0$.
\end{assumption}

\begin{theorem}\label{thm:large_n0_feasibility}
Consider the PDE system \eqref{eq:pde_main}--\eqref{eq:pde_bc} under the NN-based controller \eqref{eq:nn_controller1}--\eqref{eq:nn_controller3} constructed via the aggregate training condition \eqref{eq:training_LMI}. Assume that $p \in C^{1}([0,1])$ with $p>0$, the Sturm--Liouville eigenvalue asymptotics give $\lambda_n/n^2\to c_\lambda>0$, while Lemma~2 gives $|\beta_n|=O(\sqrt{\lambda_n})$. Assume moreover that there exist $\delta_0>0$ and $N_u\in\mathbb{N}$ such that
\begin{equation}
  -\lambda_n + q_c \le -\delta_0 < 0, \qquad \forall n \ge N_u .
  \label{eq:tail_dissipation}
\end{equation}

Fix a desired decay rate $\delta\in(0,\delta_0)$ and choose any $\gamma>0$.

The spillover matrix is $\mathcal R=\gamma\mu S_2(n_0)\tilde R_u^\top\tilde R_u$ with $S_2(n_0)=\sum_{n>n_0}\beta_n^2/\lambda_n^2$, and $S_2(n_0)\to0$ as $n_0\to\infty$ because $|\beta_n|=O(\sqrt{\lambda_n})$ (Lemma~\ref{lemma:beta_properties}).

{\color{black}
If Assumption~\ref{ass:certifiable_robust_expert} holds, then there exists an integer $N^\star\geq\max\{N_0,N_u\}$ such that, for every truncation order $n_0\geq N^\star$, one can fix $\mu\geq1/2$ satisfying $\Gamma_{n_0+1}\leq0$ so that the aggregate training condition \eqref{eq:training_LMI} is feasible.} The resulting infinite-dimensional closed-loop system is locally exponentially stable with decay rate $\delta$ in the sense of Theorem~\ref{theo}.
\end{theorem}

\begin{pf}
\noindent\textbf{Step 1 (Uniform dissipation for the truncated certificate).} 

{\color{black}Fix $\delta$ and $\gamma$. By Assumption~\ref{ass:certifiable_robust_expert}, for every $n_0\geq N_0$ there exist $P\succ0$, $\Theta\succ0$, and $\tilde N$ satisfying the containment condition and the uniform spillover-free margin \eqref{eq:uniform_margin}.} The margin quantifies the capacity of the finite-dimensional certificate to absorb spillover: for any positive semi-definite perturbation $\mathcal R\prec\eta I$,
{\color{black}
\begin{equation*}
    \tilde{R}_{\nu}^{\top} \bar{Q}_{1} \tilde{R}_{\nu}
    + \tilde{R}_{\phi}^{\top} Q_{2} \tilde{R}_{\phi}
    +\mathcal R
    \preceq -\eta I+\mathcal R\prec0.
\end{equation*}
By the change of variables and the Schur complement already used in the derivation of \eqref{eq:training_LMI}, this inequality is equivalent to $\mathcal M_H(\bar H,\bar L)\succ0$.
}

\noindent\textbf{Step 2 (Feasibility of the infinite-dimensional tail decay).} 

{\color{black}Let $n_0\geq N_u$. Before solving the training problem, fix $\mu\geq1/2$ so that $\Gamma_{n_0+1}\leq0$.} Using \eqref{eq:tail_dissipation} and $\delta<\delta_0$, we have $-\lambda_{n_0+1}+q_c+\delta\le-(\delta_0-\delta)<0$, so both conditions are met by the choice
\begin{equation}
    \mu \;=\; \max\left\{\frac{1}{2},\ \frac{\lambda_{n_0+1}}{2(\lambda_{n_0+1}-q_c-\delta)}\right\},
    \label{eq:alpha_choice}
\end{equation}
whose second argument tends to $1/2$ as $n_0\to\infty$, so $\mu$ is uniformly bounded. With this choice, $\Gamma_{n_0+1}\le0$ and $\mu\ge1/2$ guarantee $\Gamma_n\le0$ for every $n>n_0$ (see the proof of Theorem~\ref{theo}). {\color{black}Hence the residual-mode condition does not obstruct feasibility for any $n_0\geq N_u$ and requires no a posteriori optimization.}

\noindent\textbf{Step 3 ({\color{black}Feasibility of the direct-spillover training inequality}).} 

{\color{black}By Step 1, the direct-spillover component of \eqref{eq:training_LMI} holds whenever $\mathcal R=\gamma\mu S_2(n_0)\tilde R_u^\top\tilde R_u\prec\eta I$, for which it suffices that $\gamma\mu S_2(n_0)\|\tilde R_u\|^2<\eta$.} By Assumption~\ref{ass:certifiable_robust_expert}, $\|\tilde R_u\|\le\kappa_u$ uniformly in $n_0$ (by Step 2, $\mu$ is uniformly bounded) and
\begin{equation*}
    S_2(n_0)=\sum_{n>n_0}\frac{\beta_n^2}{\lambda_n^2}\ \xrightarrow[n_0\to\infty]{}\ 0,
\end{equation*}
 because $|\beta_n|=O(\sqrt{\lambda_n})$ and $\lambda_n/n^2\to c_\lambda>0$ (Lemma~\ref{lemma:beta_properties}) give $\beta_n^2/\lambda_n^2=O(n^{-2})$. Hence there exists $N^\star\ge\max\{N_0,N_u\}$ such that {\color{black}$\gamma\mu S_2(n_0)\kappa_u^2<\eta$} for all $n_0\ge N^\star$, and {\color{black}$\mathcal M_H(\bar H,\bar L)\succ0$ holds} for all such $n_0$.


{\color{black}

Therefore, there exists $N^\star\geq\max\{N_0,N_u\}$ such that, for every $n_0\geq N^\star$, the aggregate direct-spillover condition \eqref{eq:training_LMI} is feasible for the fixed choice of $\mu$ satisfying \eqref{eq:lmi2}.  Under the inverse change of variables \eqref{def:1}--\eqref{def:3}, the corresponding certificate satisfies the sufficient conditions of Theorem~\ref{theo}. Hence, the associated infinite-dimensional closed-loop system is locally exponentially stable with decay rate $\delta$.
}

This concludes the proof. \hfill $\square$
\end{pf}

\section{Numerical Implementation of the Stability-Constrained Learning Framework via ADMM}\label{section:ADMM}

This section presents the numerical implementation of the proposed stability-constrained learning framework. We formulate the joint optimization of the NN weights and stability-certificate variables and introduce an ADMM-based procedure that alternates their updates. The algorithm also verifies the stability conditions for the learned controller with its weights held fixed.

\subsection{Stability-Constrained Imitation-Learning Problem}

Given expert demonstrations consisting of state-control pairs, we train an NN controller with weights $N$ by minimizing a loss function $\mathcal{L}(N)$, which measures the discrepancy between the NN's predictions and the expert data. While common choices for $\mathcal{L}(N)$ include mean-squared error or another regression loss, the loss is generally non-convex in $N$. To balance imitation accuracy with a certified retained-state region of attraction, we formulate the following constrained optimization problem:
\begin{align}
    \minimize{N,\bar{H},\bar{L}} & \eta_1 \mathcal{L}(N) - \eta_2 \log \det(H_1)  \label{prob:learning_obj} \\
    \subjectto & \text{LMI}(\bar{H}, \bar{L}) \mathrel\succ 0, \label{prob:learning_lmi}  \\
    & \mathcal{F}(N)H = L, \label{prob:learning_consistency} \\
    & S_iH_1S_i^\top\le s_i^2,
\end{align}
\noindent where 
\begin{align*}
    H \coloneqq \begin{bmatrix}
        H_{1} & 0 \\ 0 & H_{2}
    \end{bmatrix}, \quad L \coloneqq \begin{bmatrix}
        L_{1} & L_{2} \\ L_{3} & L_{4}
    \end{bmatrix}.
\end{align*}

The objective \eqref{prob:learning_obj} balances imitation accuracy, through $\mathcal{L}(N)$, with the size of the certified retained-state ellipsoid, through $\log \det(H_1)$, where $\eta_1, \eta_2 > 0$ are weighting parameters. The parameters $\eta_{1}$ and $\eta_{2}$ are user-defined weights that directly control the trade-off between the two main objectives of the learning process: $\eta_{1}$ determines the emphasis placed on matching the expert control actions, while $\eta_{2}$ governs the emphasis on safety and robustness by prioritizing controllers with a larger certified Region of Attraction (ROA). The equality constraint $\mathcal{F}(N)H = L$ enforces consistency between the loop-transformed representation $\mathcal F(N)=\tilde N$ of the actual NN controller and the auxiliary variables used in the convex stability formulation. Since $H\succ0$, this constraint implies $\tilde N=\mathcal F(N) =LH^{-1}$.

{\color{black}
\begin{theorem}\label{theo2}
Consider the reaction--diffusion system
\eqref{eq:pde_main}--\eqref{eq:pde_ini} under the standing PDE, truncation, and local-sector assumptions of Theorem~\ref{theo}. In particular, let $n_0\geq1$ and $\delta>0$ satisfy
\begin{equation}
    -\lambda_n+q_c<-\delta,
    \qquad \forall n\geq n_0+1,
\end{equation}
and assume that the componentwise local sector bounds of Lemma~\ref{lemma:sector_condition}, computed for the corresponding NN controller, are valid on $\mathcal Z$, with $\mathcal Z\subseteq\mathcal Z_h$.

Fix $\gamma>0$ and $\mu\geq1/2$ satisfying the residual-mode condition \eqref{eq:lmi2}. Let $\tilde N=\mathcal F(N)$ denote the exact loop-transformed representation of the NN controller introduced in Lemma~\ref{lemma:loop_transformation}.

Suppose that $(N,\bar H,\bar L)$ is a feasible point of the stability-constrained learning problem \eqref{prob:learning_lmi}--\eqref{prob:learning_consistency}, with
\begin{equation}
    H_1\in\mathbb S_{++}^{n_0},
    \qquad
    H_2\in\mathbb D_{++}^{n_\phi}.
\end{equation}

Define
\begin{equation}
    P\coloneqq H_1^{-1},
    \qquad
    \Theta\coloneqq H_2^{-1}.
\end{equation}

Then $P\in\mathbb S_{++}^{n_0}$ and $\Theta\in\mathbb D_{++}^{n_\phi}$ satisfy \eqref{eq:lmi} and \eqref{eq:lmi4}. Together with the residual-mode condition \eqref{eq:lmi2}, all sufficient conditions of Theorem~\ref{theo} are satisfied.

Consequently, the origin of the complete infinite\-/dimensional closed-loop system is locally exponentially stable in $\mathcal H_{-1/2}$ with decay rate $\delta$. Moreover, $\mathcal E(P,\gamma)$ defined in \eqref{eq:set_invariance} is positively invariant and constitutes an inner approximation of the region of attraction. In particular, there exists $M\geq1$, independent of the initial condition, such that
\begin{equation}
    \|z(t,\cdot)\|_\star
    \leq
    M e^{-\delta t}\|z(0,\cdot)\|_\star,
    \qquad \forall t\geq0,
\end{equation}
for every $z(0,\cdot)\in\mathcal E(P,\gamma)$.
\end{theorem}

\begin{pf}
Let $(N,\bar H,\bar L)$ be a feasible point of the stability-constrained learning problem. Since $H_1\in\mathbb S_{++}^{n_0}$ and $H_2\in\mathbb D_{++}^{n_\phi}$, both matrices are nonsingular. Hence, the recovered matrices
\begin{equation}
    P\coloneqq H_1^{-1},
    \qquad
    \Theta\coloneqq H_2^{-1}
\end{equation}
satisfy $P\in\mathbb S_{++}^{n_0}$ and $\Theta\in\mathbb D_{++}^{n_\phi}$.

By feasibility of the consistency constraint and the identity $\tilde N=\mathcal F(N)$, we have
\begin{equation}
    \tilde N H=L.
\end{equation}

Using the block structures of $\tilde N$, $H$, and $L$, this identity implies precisely the change-of-variable relations \eqref{def:2}--\eqref{def:3}. Therefore, the auxiliary variables $\bar L$ appearing in \eqref{eq:training_LMI} are consistent with the exact loop-transformed representation of the NN controller.

Consider now the stability component of \eqref{eq:training_LMI}. By the Schur-complement equivalence \eqref{eq:schur_spillover} and the invertibility of $H_1$ and $H_2$, the congruence transformation used in deriving $M_H$ can be reversed. Using $P=H_1^{-1}$, $\Theta=H_2^{-1}$, and the change-of-variable relations \eqref{def:2}--\eqref{def:3}, this recovers the Schur-complement representation of the direct-spillover stability condition \eqref{eq:lmi}. Since $\Theta^{-1}=H_2\succ0$, the corresponding Schur complement is reversible. Hence \eqref{eq:lmi} holds for the exact loop-transformed NN controller.

Similarly, the containment component of \eqref{eq:training_LMI} satisfies, for every $i\in\{1,\ldots,n_z\}$,
\begin{equation}
    \begin{bmatrix}
        s_i^2 & S_i^\top H_1\\
        H_1S_i & H_1
    \end{bmatrix}
    \succeq0.
\end{equation}

Since $P=H_1^{-1}$, applying the nonsingular congruence transformation $\diag(1,P)$ gives
\begin{equation}
    \begin{bmatrix}
        s_i^2 & S_i^\top\\
        S_i & P
    \end{bmatrix}
    \succeq0,
\end{equation}
which is precisely the containment condition \eqref{eq:lmi4}.

Finally, the residual-mode condition \eqref{eq:lmi2} holds by hypothesis. Hence \eqref{eq:lmi}, \eqref{eq:lmi2}, and \eqref{eq:lmi4} are all satisfied by the recovered matrices $P$ and $\Theta$ and the exact loop-transformed representation of the NN controller.

Therefore, all sufficient conditions of Theorem~\ref{theo} hold. The positive invariance of $\mathcal E(P,\gamma)$ and the local exponential stability of the complete infinite-dimensional closed-loop system in $\mathcal H_{-1/2}$ with decay rate $\delta$ follow directly from Theorem~\ref{theo}.

This concludes the proof. \hfill $\square$
\end{pf}
}


\subsection{ADMM-Based Training Algorithm}

We employ the Alternating Direction Method of Multipliers (ADMM) \cite{yin2021imitation} to decouple the non-convex learning problem from the convex constraints. The augmented Lagrangian is defined as:
\begin{multline*}
    \mathcal{L}_a(N, \bar{H}, \bar{L}, Y) = \eta_1 \mathcal{L}(N) - \eta_2 \log \det(H_1) \\ + \operatorname{tr}\left(Y^{\top} (\mathcal{F}(N)H - L)\right) + \frac{\rho}{2} \left\|\mathcal{F}(N)H - L\right\|_F^2,
\end{multline*}
\noindent where $Y$ is the Lagrange multiplier, $\rho > 0$ is a regularization parameter, and $\|\cdot\|_F$ denotes the Frobenius norm. In contrast, $\rho$ is not related to the control objectives but is a regularization parameter for the ADMM optimization algorithm. It functions as a penalty term that helps managing the constraints and influences the numerical convergence of the solver. The ADMM iterations proceed as shown in Algorithm \ref{alg:safe_il_pde}.

\begin{algorithm*}
\caption{Safe Imitation Learning for PDE Systems}
\label{alg:safe_il_pde}
\renewcommand{\algorithmicrequire}{\textbf{Input:}}
\renewcommand{\algorithmicensure}{\textbf{Output:}}
\renewcommand{\algorithmiccomment}[1]{\hfill #1}
\algnewcommand\algorithmicinitialize{\textbf{Initialize:}}
\algnewcommand\Initialize{\item[\algorithmicinitialize]}

\begin{algorithmic}[1]
\Require Expert data pairs, weighting parameters $\eta_1,\eta_2$, penalty $\rho>0$, training epochs $E$, and fixed parameters $\gamma>0$ and $\mu\geq1/2$ satisfying \eqref{eq:lmi2}.
\Initialize NN weights $N^0$, LMI variables $(\bar H^0,\bar L^0)$, Lagrange multiplier $Y^0$, and iteration counter $k\gets0$.

\Repeat
    \State \Comment{\textit{Step 1: NN Parameter Update}}
    \State Set $N^{k+1}\gets N^k$.
    \For{$e=1\to E$} \Comment{Inner training loop}
        \State Update the NN weights by performing one gradient step
        on the augmented loss:
        \State $N^{k+1}\gets\text{AdamOptimizer} \big( \nabla_N \mathcal L_a( N^{k+1},\bar H^k, \bar L^k,Y^k) \big)$.
    \EndFor

    \State \Comment{\textit{Step 2: Stability Certificate Update}}
    \State Solve the convex optimization problem for the stability variables:
    \State $(\bar H,\bar L)^{k+1} \gets \arg\min_{\bar H,\bar L} \mathcal L_a( N^{k+1}, \bar H, \bar L, Y^k)$
    \Statex $\hspace{\algorithmicindent}\text{subject to}\quad\operatorname{LMI}(\bar H,\bar L)\succ0$.

    \State \Comment{\textit{Step 3: Lagrange Multiplier Update}}
    \State $Y^{k+1} \gets Y^k+ \rho\big(\mathcal F(N^{k+1})H^{k+1}-L^{k+1}\big)$.

    \State \Comment{\textit{Step 4: Fixed-Controller Certificate Verification}}
    \State With $N^{k+1}$ fixed, recompute its exact loop-transformed representation $\mathcal F(N^{k+1})$ and solve the certificate problem
    \State $(\bar H^\star,\bar L^\star) \gets \arg\min_{\bar H,\bar L}-\log\det(H_1)$
    \Statex $\hspace{\algorithmicindent}\text{subject to}\quad \operatorname{LMI}(\bar H,\bar L)\succ0,\qquad \mathcal F(N^{k+1})H=L$.
    \State $k\gets k+1$.
\Until{the ADMM convergence criterion is satisfied and the fixed-controller certificate problem is feasible.}

\State \Return Converged NN weights $N^\star$ and the exact stability certificate $(\bar H^\star,\bar L^\star)$.

\Statex \vspace*{-1ex}
\end{algorithmic}
\end{algorithm*}

Because the learning problem is non-convex, convergence of the ADMM iterations does not by itself constitute a stability certificate. A trained controller is declared certified only if the final, fixed NN is used to recompute its actual loop-transformed matrix $\tilde N^\star=\mathcal F(N^\star)$, its local sector bounds,  and all matrices appearing in \eqref{eq:lmi}--\eqref{eq:lmi4}, and if the resulting finite- and infinite-dimensional verification inequalities are satisfied within the prescribed numerical tolerances.


\section{Results} \label{sec:results}

To illustrate the results, we consider a PDE system \eqref{eq:pde_main}--\eqref{eq:pde_bc} subject to the homogeneous transformations with the case where the diffusion coefficient is $p(x)=1$. The spectrum of the operator $\mathcal{A}$ can then be determined analytically, and it is well-documented (see, e.g., \cite{boyce2021elementary}). The eigenvalues are given by $\lambda_{n} = \left(n - \frac{1}{2}\right)^2 \pi^2$, with corresponding eigenfunctions $\Phi_{n}(x) = \sqrt{2} \cos\left(\left(n - \frac{1}{2}\right)\pi x\right)$, $\forall n\geq1$. For this example, $p(1)=1$, and Lemma~\ref{lemma:beta_properties} gives $|\beta_n|=\sqrt{2\lambda_n}$. Therefore, $S_2(n_0)=2\sum_{n=n_0+1}^{\infty}\frac{1}{\lambda_n}$. This provides an exact expression for the spillover coefficient used in the numerical verification. A finite-dimensional model suitable for controller design is then derived by truncating this system to its first $n_{0}$ modes, as described by \eqref{eq:sys_trunc}.

All numerical experiments were performed in \texttt{Python 3.9} on a computer equipped with an \texttt{Apple M4 Pro} processor and \texttt{24}~GB of unified memory. The NN was trained using \texttt{TensorFlow 2.12.0}, the semidefinite programs were solved using \texttt{CVXPY} with the \texttt{MOSEK 11} solver, and the MPC optimization problems were formulated as convex quadratic programs (QPs) and solved using \texttt{CVXOPT 1.3.2}.

Using this model, the expert controller is defined as the following linear MPC:
\begin{align*}
   \hspace{-0.85cm} \displaystyle\min_{\{u_k\}_{k=0}^{N_p-1}}\;\; & \sum_{k=0}^{N_{p}-1} \left( Z_{n_{0},k}^{\top}\bar{Q}Z_{n_{0},k} + u_{k}^{\top}\bar{R}u_{k} \right) + Z_{n_{0},N_{p}}^{\top}\bar{P}Z_{n_{0},N_{p}} \\
   \hspace{-0.85cm} \mbox{s.t.} \;\;  & Z_{n_{0},k+1} = A_{d}Z_{n_{0},k} + B_{d}u_{k}, \\
               & Z_{n_{0},0} = Z_{0}, \\
               & u_{k} \in \mathcal{U},\, \forall k \in [0, \dots, N_{p}-1], \\
               & Z_{n_{0},k} \in \mathcal{Z},\, \forall k \in [0, \dots, N_{p}-1], \\
               & Z_{n_{0},N_{p}} \in \mathcal{Z}_{t}, 
\end{align*}
\noindent where $N_{p}$ is the prediction horizon, $\bar P$ is the terminal cost matrix, and $\mathcal{Z}_{t}$ is the terminal constraint set. The terminal feedback gain $K_{\mathrm{lqr}}$ is obtained from the discrete-time algebraic Riccati equation associated with the weighting matrices $\bar{Q} = 2I_{n_{0}}$ and $\bar{R} = 1$. The terminal set is chosen as $\mathcal{Z}_{t} = \mathcal Z_\infty$, where $\mathcal Z_\infty$ is the maximal positively invariant (MPI) set of the closed-loop system $Z_{n_0,k+1}=(A_d-B_dK_{\mathrm{lqr}})Z_{n_0,k}$, that is
\begin{equation*}
    \mathcal Z_\infty = \left\{ Z_{n_0}\in\mathcal Z:\; (A_{d} - B_{d}K_{\mathrm{lqr}})Z_{n_0}\in\mathcal Z_\infty \right\},
\end{equation*}
\noindent which is computed offline by fixed-point iteration.

\subsection{NN Controller for the Unstable Reaction-Diffusion PDE}

We start with a more conservative approach. Setting $q_{c}=22.5$ and $\delta=0.01$, the tail dissipation relation $-\lambda_{n}+q_{c}<-\delta<0$ for all $n\geq n_{0}+1$ yields the truncation order $n_0=2$, corresponding to the two unstable modes. For the MPC controller, we impose the state constraint set
\begin{equation*}
    \mathcal{Z} = \left\{ z\in\mathbb{R}^{2} ~\middle|~ |z_{1}|\leq2,\, |z_{2}|\leq30 \right\},
\end{equation*}
and the control constraint set
\begin{equation*}
    \mathcal{U} = \left\{ u\in\mathbb{R} ~\middle|~ |u|\leq30 \right\}.
\end{equation*}

Unless otherwise stated, all learning experiments in this subsection use $\rho=1$, $\eta_{1}=1$, and the prescribed decay rate $\delta=0.01$, while the MPC uses a sampling period of $T_s=0.02~\mathrm{s}$ with zero-order-hold (ZOH) actuation.

The following experiments investigate three different aspects of the proposed formulation. First, the weighting parameter $\eta_2$ is varied to illustrate the trade-off between imitation accuracy and the size of the certified region of attraction. Second, the tail weight $\gamma$ is varied while keeping the trained NN controller fixed, thereby isolating its effect on the infinite-dimensional stability certificate. Finally, the truncation order is increased from $n_0=2$ to $n_0+1$ by including the first stable mode in the retained model, allowing the effect of the residual-tail approximation on the conservatism of the certificate to be assessed.

\paragraph{Influence of the ROA weight $\eta_2$:}

\begin{figure}
    \centering
    \includegraphics[width=\linewidth]{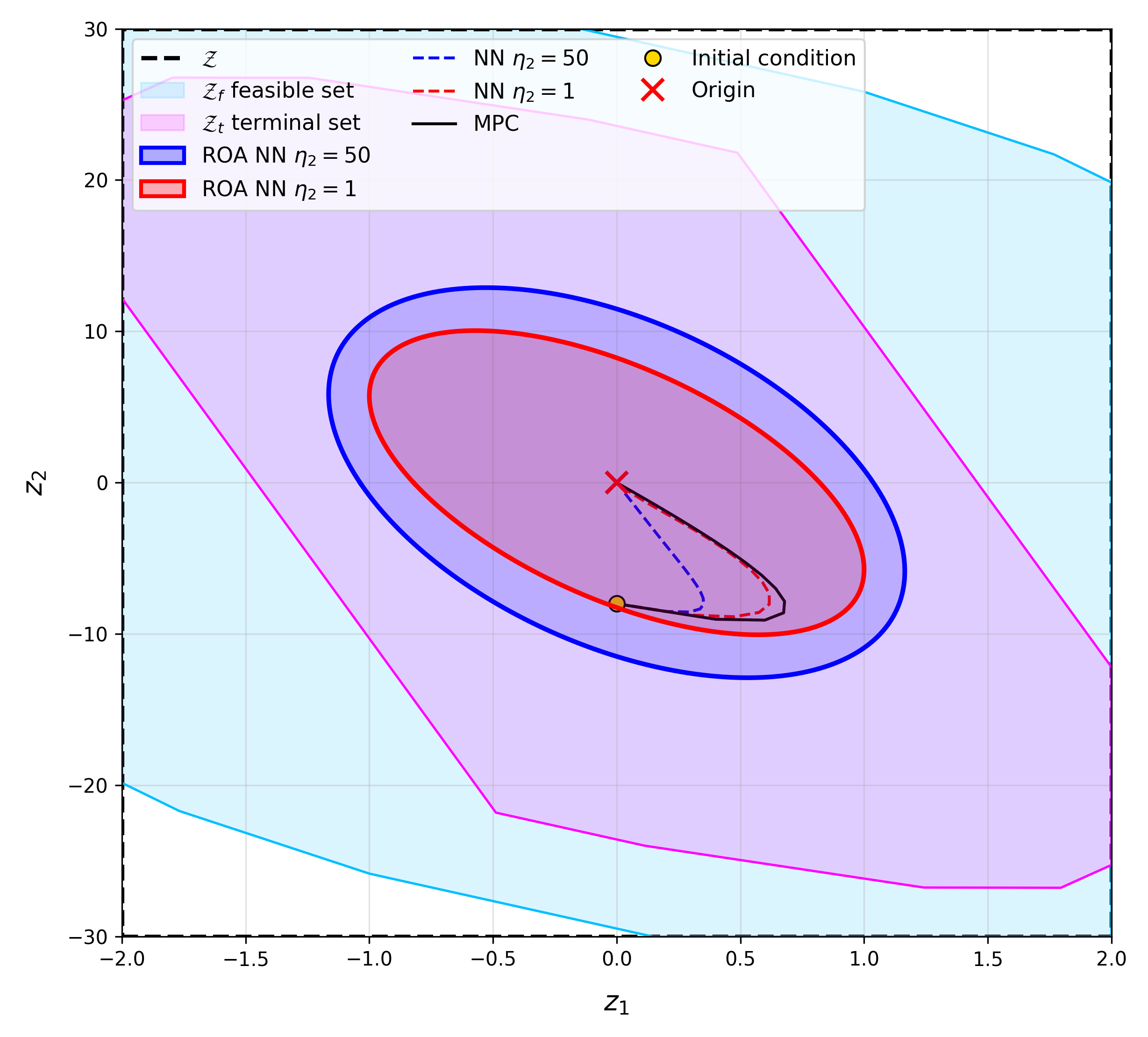}
    \caption{Truncated-state trajectories and certified regions of attraction obtained for $\eta_2=1$ and $\eta_2=50$.}
    \label{fig:eta2_comparison}
\end{figure}

Figure~\ref{fig:eta2_comparison} compares the closed-loop responses and certified ellipsoids obtained for different values of $\eta_2$. Recall that the stability-constrained learning problem minimizes \eqref{prob:learning_obj} with $P=H_1^{-1}$. Hence, increasing $\eta_2$ increases the relative importance assigned to enlarging the certified retained-state ellipsoid, whereas decreasing $\eta_2$ places relatively greater emphasis on reproducing the expert policy.

This trade-off can be clearly observed in Figure~\ref{fig:eta2_comparison}. For $\eta_2=50$, the optimization places a substantial weight on the ellipsoidal-volume term and returns a larger certified ROA. This increase in the certified region is obtained at the expense of a larger deviation from the expert policy, and the resulting closed-loop state trajectories exhibit a slower transient response. In contrast, for $\eta_2=1$, imitation accuracy dominates the optimization. The resulting NN controller remains closer to the MPC expert and consequently reproduces its faster transient behavior more accurately, although the certified ROA is smaller.

\begin{figure}
    \centering
    \includegraphics[width=\linewidth]{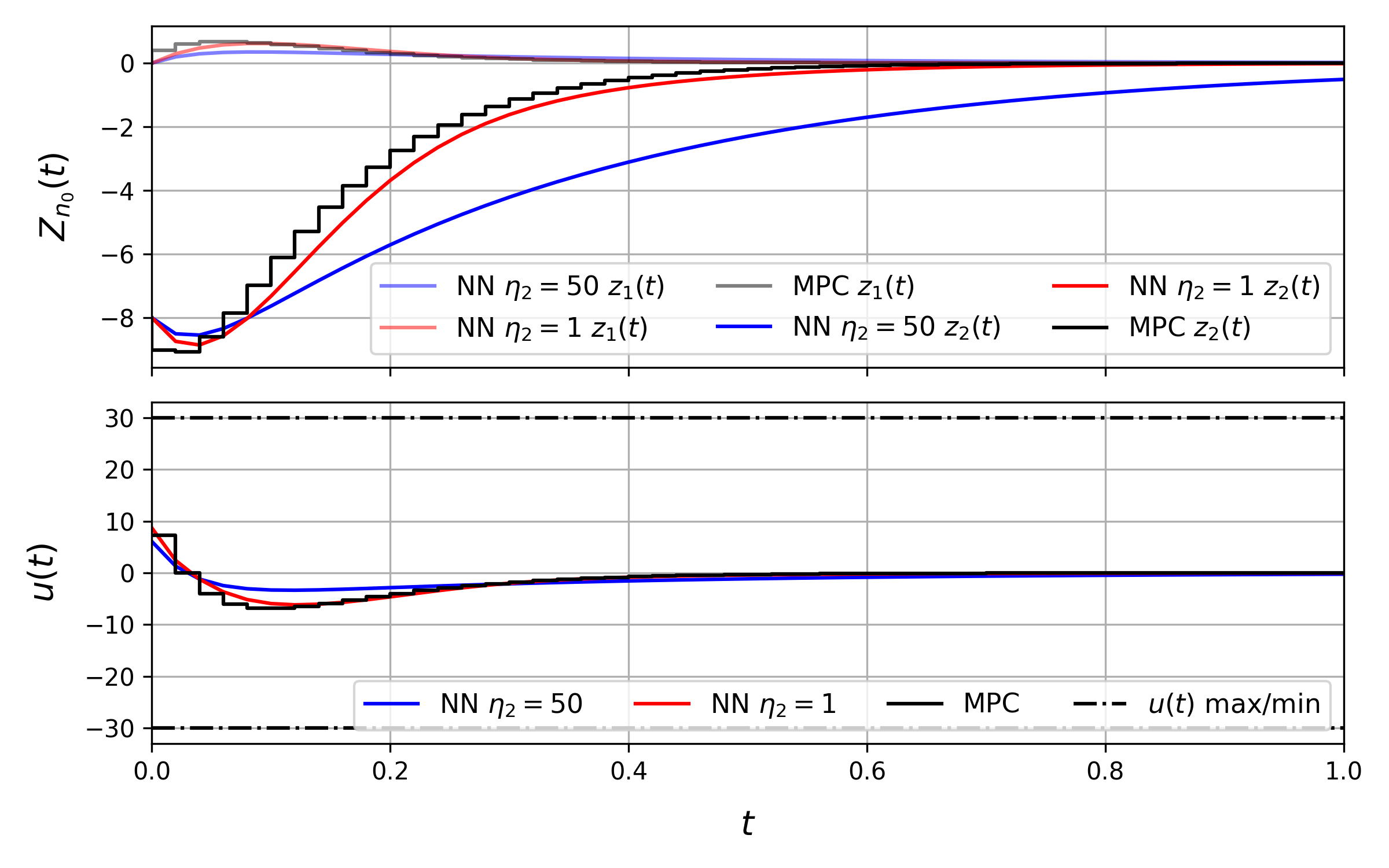}
    \caption{Time-domain comparison of the truncated-state states and boundary control input for the NN controllers trained with $\eta_2=1$ and $\eta_2=50$  together with the MPC expert.}
    \label{fig:eta2_comparison_control}
\end{figure}


Therefore, a larger value of $\eta_2$ permits the optimizer to accept some degradation in imitation accuracy if doing so allows a larger set of states to satisfy the stability and local-sector constraints. Conversely, a smaller $\eta_2$ provides less incentive to enlarge the certified set, allowing the optimization to remain closer to the MPC demonstrations.


Figure~\ref{fig:eta2_comparison_control} complements this comparison by presenting the retained-state trajectories and boundary-control inputs generated by the NN controllers trained, together with those generated by the MPC expert.

\begin{figure}
    \centering
    \includegraphics[width=\linewidth]{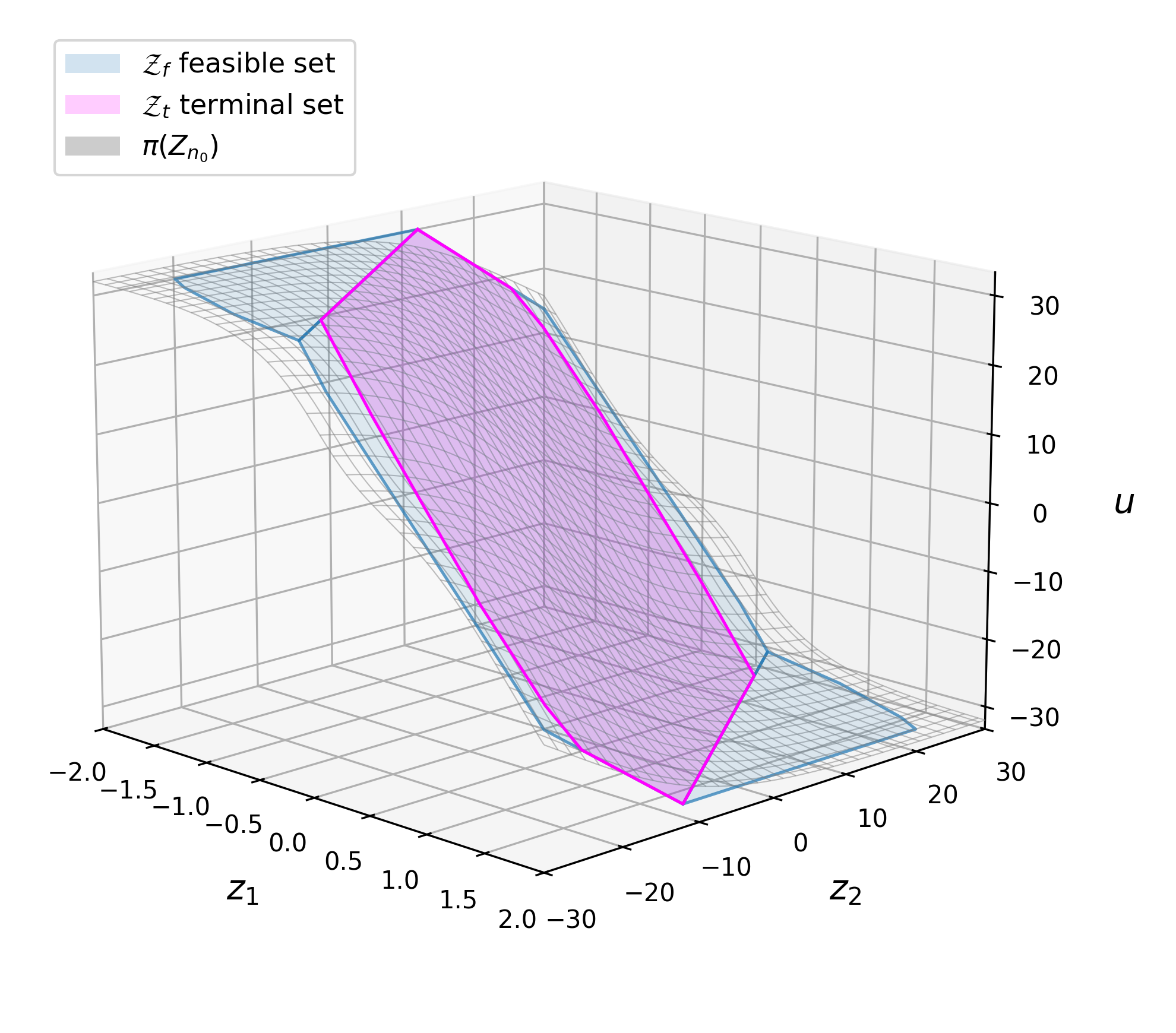}
    \caption{Learned NN policy surface $u=\pi(Z_{n_0})$ over the admissible state domain.}
    \label{fig:policy_surface}
\end{figure}

Figure~\ref{fig:policy_surface} complements the trajectory-based comparison by visualizing the learned feedback map over the retained-state domain. While the closed-loop trajectories show the NN behavior only along particular state evolutions, the policy surface displays the control action generated throughout the admissible region, thereby the Figure illustrates how the learned NN maps each retained-state pair $(Z_{n_0})$ to the corresponding control action.

\paragraph{Influence of the residual-tail weight $\gamma$:}

\begin{figure}
    \centering
    \includegraphics[width=\linewidth]{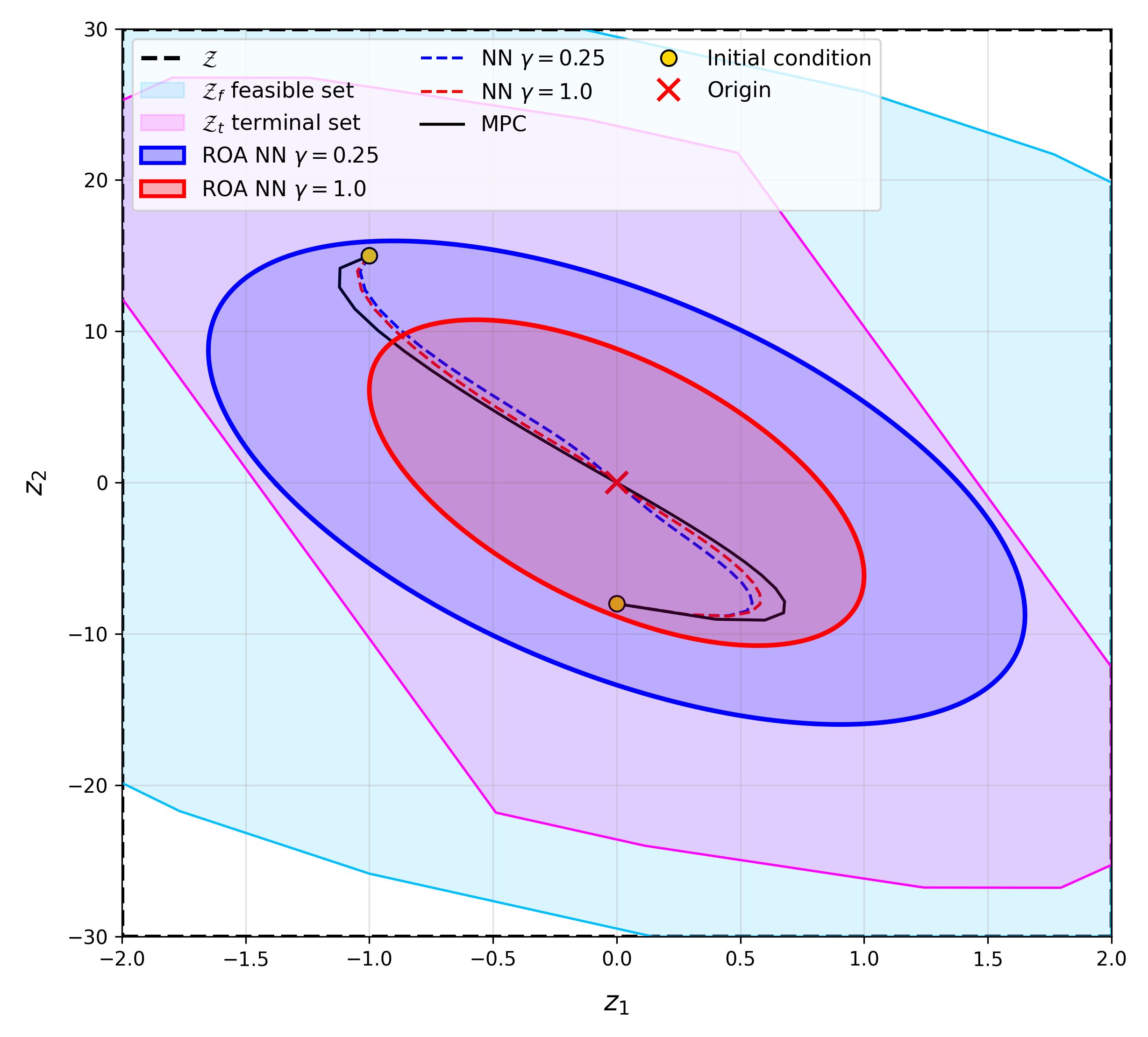}
    \caption{Truncated-state trajectories and certified regions of attraction for $\gamma=0.25$ and $\gamma=1$.}
    \label{fig:gamma_comparison}
\end{figure}

Figure~\ref{fig:gamma_comparison} investigates a conceptually different effect by comparing $\gamma=1$ and $\gamma=0.25$. Consequently, $\gamma$ modifies the stability certificate.

The effect of $\gamma$ follows directly from the composite Lyapunov functional \eqref{eq:V1}--\eqref{eq:V2}. The scalar $\gamma$ determines the relative weight assigned to the residual modal energy. More importantly for the LMI condition, it appears explicitly in the direct-spillover term \eqref{eq:spill}

Thus, for fixed $\mu$, $n_0$, and controller weights, increasing $\gamma$ increases the positive-semidefinite perturbation that must be absorbed by the finite-dimensional dissipation margin.

This effect is observed in Figure~\ref{fig:gamma_comparison}. With $\gamma=1$, a larger contribution from the residual dynamics is included in the certificate, resulting in a more restrictive stability inequality and, consequently, a smaller certified ROA. Reducing the residual weight to $\gamma=0.25$ decreases the spillover contribution entering the LMI and allows a larger retained-state ellipsoid to be certified. 

\begin{figure}
    \centering
    \includegraphics[width=\linewidth]{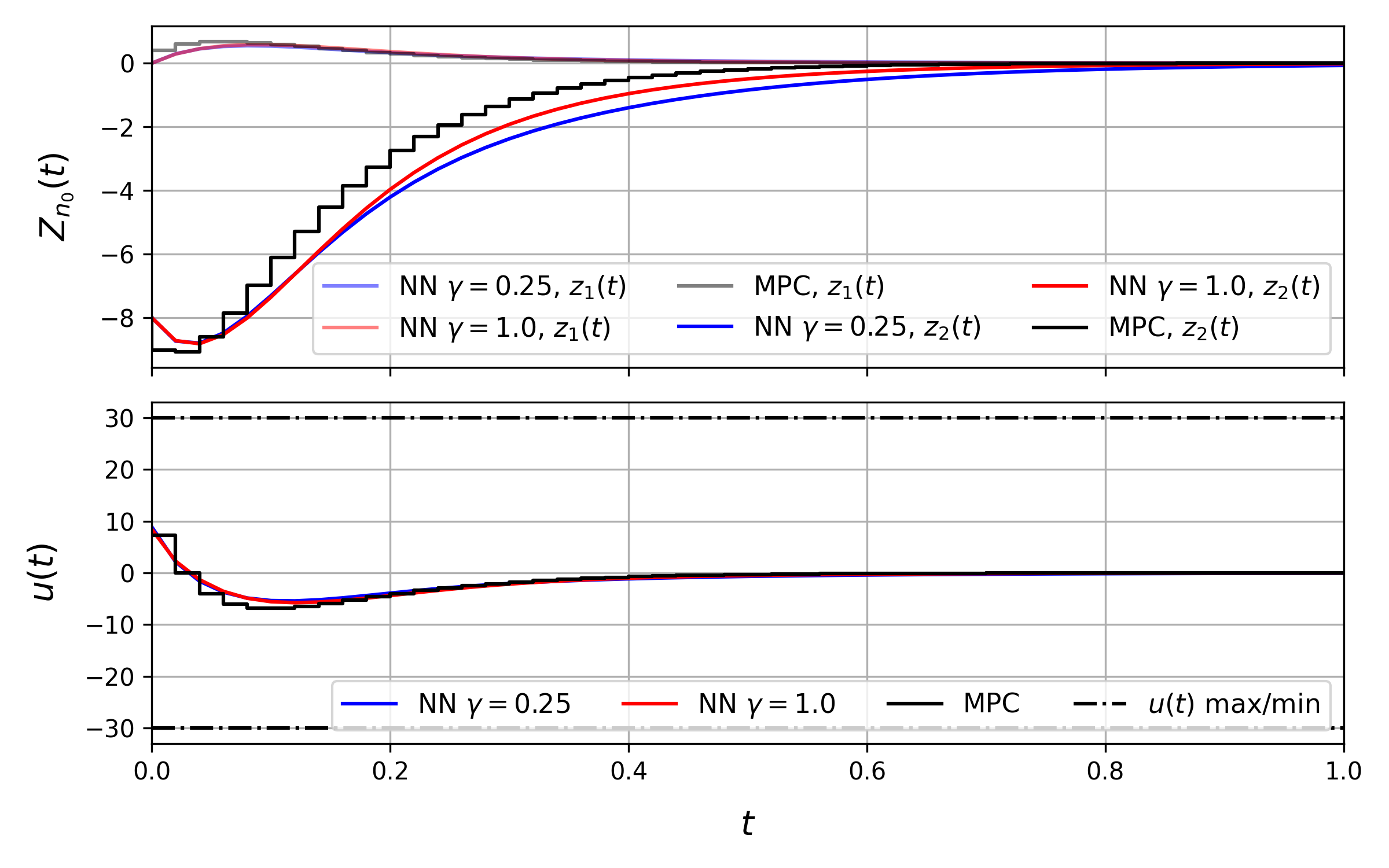}
    \caption{Time-domain truncated-state and boundary-control responses for $\gamma=0.25$ and $\gamma=1$, and compared with the MPC expert.}
    \label{fig:gamma_comparison_control}
\end{figure}

The comparison also illustrates an important distinction between closed-loop performance and certification. Changing $\gamma$ does not itself make the controller faster or slower. Instead, it changes how strongly the residual infinite-dimensional state is weighted in the Lyapunov certificate and how much robustness margin must be allocated to account for spillover. Hence, substantially different certified regions can be obtained for essentially the same state and control trajectories.


Figure~\ref{fig:gamma_comparison_control} presents the retained-state trajectories and boundary-control input generated by the fixed NN controller for the two certification values, together with the response of the MPC expert. Since the same NN controller is used in both cases, the corresponding NN trajectories and control inputs coincide, changing \(\gamma\) affects only the stability certificate.

\paragraph{Influence of the truncation order and inclusion of a stable mode:}

The third experiment investigates the effect of increasing the truncation order. Starting from $n_0=2$, for which only the two unstable modes are retained, we include the next mode and set $n_0+1$. The third mode is already open-loop stable and therefore does not need to be included merely to stabilize the truncated system. Its inclusion is instead motivated by the reduction of the uncertainty associated with the neglected infinite-dimensional tail. Therefore, for the $n_0=3$ case, we define the admissible domain of the third retained mode as $\mathcal{S} = \{z_3\in\mathbb{R}\mid |z_3|\leq30\}$, and construct the augmented state constraint set $\mathcal{Z}^{(3)}\coloneqq\mathcal{Z}\times\mathcal{S} \subset\mathbb{R}^3$, with associated terminal set $\mathcal{Z}^{(3)}_f\subset\mathbb{R}^3$.

The residual spillover coefficient $S_2(n_0)$ decreases as the truncation order is increased, since fewer modes remain in the unmodeled tail. Therefore, when the third mode is included in the retained model, its contribution is removed from the residual spillover bound. As a result, the direct-spillover term entering the stability LMI becomes less restrictive, although the finite-dimensional controller representation is enlarged to account for the additional retained state.

\begin{figure}
    \centering
    \includegraphics[width=\linewidth]{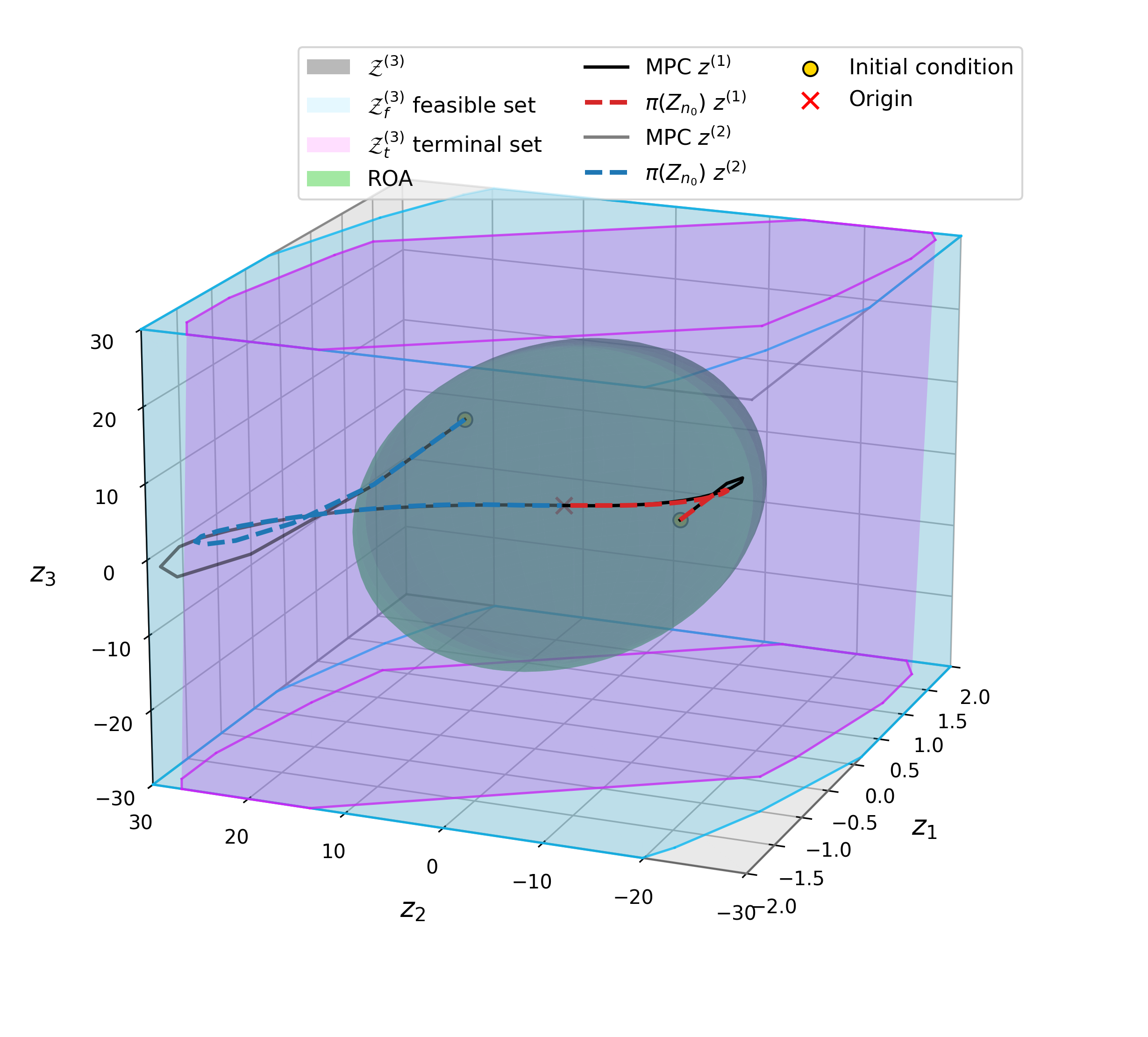}
    \caption{Three-dimensional retained-state trajectories and certified region of attraction for $n_0+1$ and $\gamma=1$.}
    \label{fig:n0_comparison}
\end{figure}
\begin{figure}
    \centering
    \includegraphics[width=\linewidth]{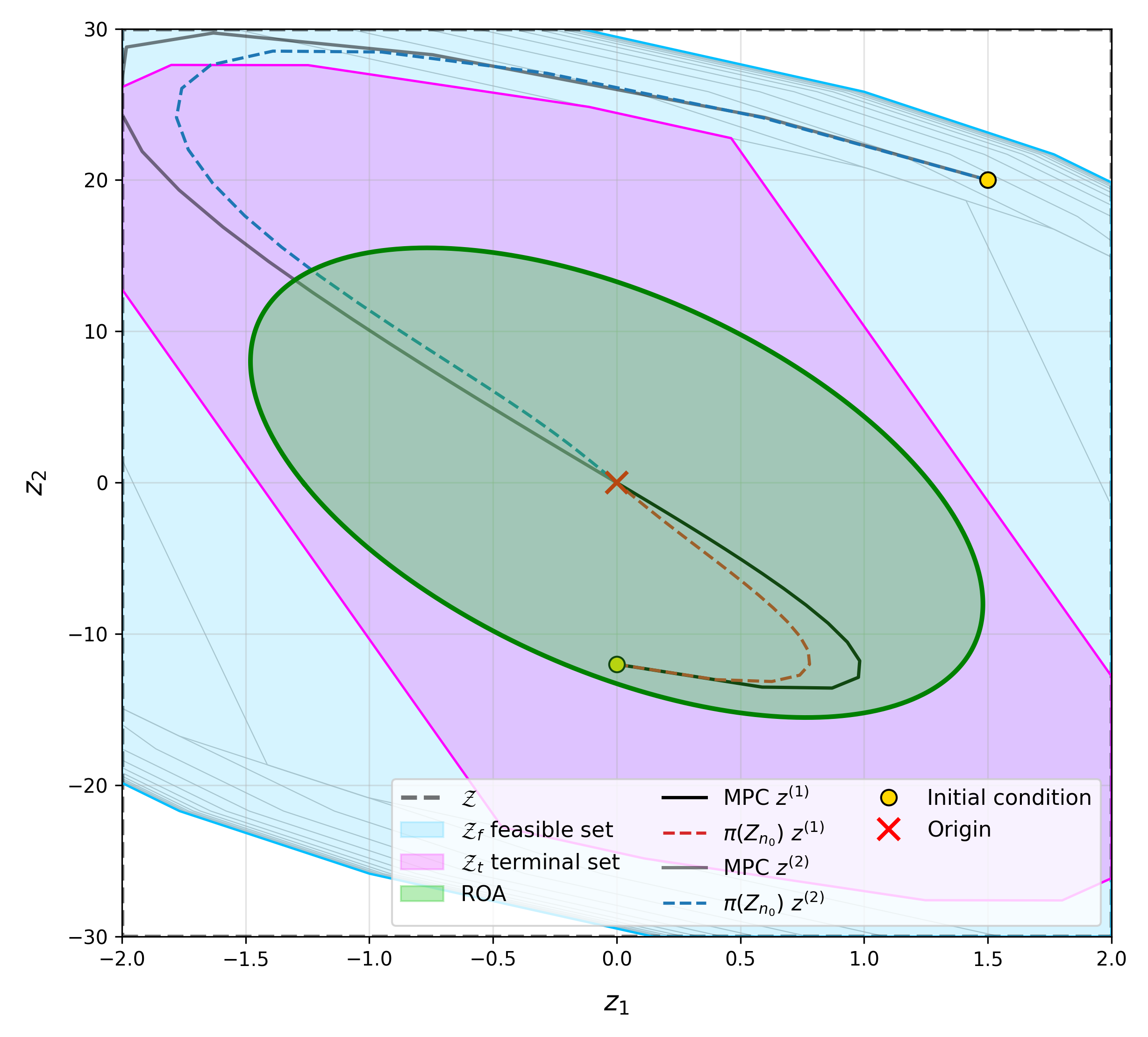}
    \caption{Projection of the $n_0+1$ certified region of attraction and closed-loop trajectories onto the $(z_1,z_2)$ plane.}
    \label{fig:n0_comparison_2}
\end{figure}

Figures~\ref{fig:n0_comparison} and~\ref{fig:n0_comparison_2} show the certified region obtained by retaining the first stable mode, i.e., by increasing the truncation order from $n_0=2$ to $n_0=3$. A direct comparison with Figure~\ref{fig:eta2_comparison} and Figure~\ref{fig:gamma_comparison}, where $n_0=2$, shows that, for the same value $\gamma=1$, the projection of the $n_0=3$ certified region onto the common $(z_1,z_2)$-plane has a larger area than the certified region obtained with $n_0=2$. The result illustrates that including additional stable modes can reduce the conservatism associated with the neglected dynamics. Although the finite-dimensional optimization problem becomes larger, a larger fraction of the plant dynamics is represented explicitly and a smaller residual contribution must be absorbed through the spillover robustness margin.

There is also a complementary interpretation in terms of the residual-mode condition. When $n_0$ is increased, the first neglected eigenvalue changes from $\lambda_3$ to $\lambda_4$. Since the Sturm--Liouville eigenvalues increase with the modal index, the first residual mode is then more strongly dissipative. The condition \eqref{eq:Gamma} therefore becomes easier to satisfy as the first neglected mode moves further into the stable part of the spectrum. If $\mu$ is recomputed for each truncation order, this effect can further reduce the conservatism of the residual-mode estimate.

The result in Figure~\ref{fig:n0_comparison} therefore highlights a fundamental trade-off in the proposed approach. A small truncation order yields a lower-dimensional NN input and a smaller finite-dimensional certificate, but leaves a larger portion of the plant dynamics in the residual subsystem. The resulting value of $S_2(n_0)$ is larger, and a greater stability margin is required to account for spillover. Increasing $n_0$ raises the dimension of the retained model and of the learning problem, but reduces $S_2(n_0)$ and transfers part of the previously uncertain residual dynamics into the explicitly certified finite-dimensional subsystem.

Table~\ref{tab:controller_computation_time} compares the online computation time per control action for models retaining $n_0$ and $n_0+1$ modes. For the NN, this time corresponds to one forward evaluation of the policy, for the MPC, it includes updating the state-dependent quadratic program terms, solving the finite-horizon optimization problem. Using 100 feasible initial conditions and 51 evaluations per trajectory, the neural policies required, on average, $8.33$ and $8.75~\mu\mathrm{s}$, compared with $17.276$ and $17.191~\mathrm{ms}$ for the MPC, respectively. The ratios of mean computation times are approximately $2{,}074$ and $1{,}966$, demonstrating a substantial reduction in online computational cost for the tested implementations. Offline controller construction and plant simulation are excluded from these measurements. Although the mean MPC computation times are below the $20~\mathrm{ms}$ sampling period, these trajectory-averaged measurements do not establish a worst-case execution-time guarantee.

\begin{table}[htbp]
    \centering
    \caption{Online computation time per control action.
    Values are the mean $\pm$ sample standard deviation
    across 100 trajectory-averaged times, with 51 evaluations
    per trajectory. The last column is the ratio of the
    mean MPC time to the mean NN time.}
    \label{tab:controller_computation_time}
    \begin{tabular}{lccc}
        \toprule
        Retained modes
        & NN ($\mu\mathrm{s}$)
        & MPC ($\mathrm{ms}$)
        & $\overline{t}_{\mathrm{MPC}}/
           \overline{t}_{\mathrm{NN}}$ \\
        \midrule
        $n_0$
        & $8.33 \pm 1.74$
        & $17.276 \pm 0.520$
        & $2074$ \\
        $n_0+1$
        & $8.75 \pm 3.66$
        & $17.191 \pm 0.523$
        & $1966$ \\
        \bottomrule
    \end{tabular}
\end{table}


\section{Conclusion} \label{sec:conclusion}

In this paper, we introduced a constrained imitation-learning framework for the design of NN boundary controllers for reaction--diffusion PDEs. The proposed method combines a finite-dimensional spectral approximation, an expert MPC policy, and quadratic constraints describing the NN nonlinearities. A composite Lyapunov functional explicitly accounts for the residual modes and provides a regional exponential-stability certificate in a negative-order Hilbert norm associated with the Sturm--Liouville operator. For the certified configurations, the numerical results demonstrate that the trained NN reproduces the expert policy with substantially lower online computational cost while satisfying the proposed stability conditions.

From the results of this work, new questions arise that open paths for new research. A primary objective is to determine where the framework can be extended to tackle more complex PDE systems with disturbances and uncertainties, or even observers. Another promising direction is to move beyond IL by combining these stability guarantees with reinforcement learning, enabling the controller to discover policies that meet certified regional stability.

\printcredits

\appendix

\section{Well-posedness of the closed-loop system}
\label{app:well_posedness}

This appendix establishes the well-posedness result used in Theorems~\ref{theo} and~\ref{theo2}.

The argument is formulated in the modal state space $\mathcal H_{-1/2}$ introduced in Remark~\ref{rem:star_norm}. It is consistent with the negative-order state-space approach used in \cite[Appendix]{Lhachemi26}, but is explicitly presented here because the boundary conditions and feedback law are different.

\begin{lemma}[Closed-loop well-posedness]
\label{lem:closed_loop_wellposedness}
Consider the modal system
\begin{equation}
    \dot z_n(t) =(-\lambda_n+q_c)z_n(t)+\beta_n\pi\bigl(Z_{n_0}(t)\bigr),\qquad n\geq 1,\label{eq:appendix_closed_loop}
\end{equation}
where $Z_{n_0}(t)=\begin{bmatrix}z_1(t)&\cdots&z_{n_0}(t)\end{bmatrix}^{\top}$.

Assume that $\pi:\mathbb R^{n_0}\to\mathbb R$ is globally Lipschitz and globally bounded, and that
\begin{equation}
    \sum_{n=1}^{\infty}\frac{\beta_n^2}{\lambda_n^2}<\infty.
    \label{eq:appendix_beta_summability}
\end{equation}

Then, for every initial condition $z^0\in\mathcal H_{-1/2}$, system \eqref{eq:appendix_closed_loop} admits a unique global mild solution $z\in C\bigl([0,\infty);\mathcal H_{-1/2}\bigr)$.

Moreover, $z(t)\in L^2(0,1)$ for every $t>0$, and the solution depends continuously on its initial condition on every bounded time interval.
\end{lemma}
\begin{pf}
The dynamics of the first $n_0$ modes can be written as
\begin{equation}
    \dot{Z}_{n_0}(t) = A_{n_0}Z_{n_0}(t) + B_{n_0}\pi\bigl(Z_{n_0}(t)\bigr).
    \label{eq:appendix_retained_system}
\end{equation}

Since $\pi$ is globally Lipschitz, the vector field in \eqref{eq:appendix_retained_system} is globally Lipschitz. Consequently, for every $Z_{n_0}(0)\in\mathbb R^{n_0}$, there exists a unique solution $Z_{n_0}\in C^{1}([0,\infty);\mathbb R^{n_0})$.

Define $ u(t)=\pi\bigl(Z_{n_0}(t)\bigr)$. By the global boundedness of $\pi$, there exists $M_\pi>0$ such that $|u(t)|\leq M_\pi$,  $t\geq 0.$ For every $n\geq1$, define
\begin{equation}
    z_n(t) = e^{-(\lambda_n-q_c)t}z_n^0 +
    \beta_n\int_0^t e^{-(\lambda_n-q_c)(t-s)}u(s)\,ds.
    \label{eq:appendix_mild_solution}
\end{equation}

Since $\lambda_n\to+\infty$, there exists $N\in\mathbb N$ such that $\lambda_n-q_c\geq\frac{\lambda_n}{2}$, $n\geq N$. Hence, for $n\geq N$,
\begin{equation}
    \left|\beta_n\int_0^te^{-(\lambda_n-q_c)(t-s)}u(s)\,ds\right| \leq
    \frac{2M_\pi|\beta_n|}{\lambda_n}.
    \label{eq:appendix_input_bound}
\end{equation}

Consequently,
\begin{equation}
    \sum_{n=N}^{\infty}\frac{1}{\lambda_n}\left|\beta_n\int_0^te^{-(\lambda_n-q_c)(t-s)}u(s)\,ds\right|^2\leq4M_\pi^2\sum_{n=N}^{\infty}\frac{\beta_n^2}{\lambda_n^3}<\infty,
\end{equation}
where the convergence follows from \eqref{eq:appendix_beta_summability}. Together with $z^0\in\mathcal H_{-1/2}$ and the finitely many modes $n<N$, this proves that $z(t)\in\mathcal H_{-1/2}$.

The continuity of $t\mapsto z(t)$ in $\mathcal H_{-1/2}$ follows from \eqref{eq:appendix_mild_solution} and dominated convergence. Uniqueness follows from the uniqueness of \eqref{eq:appendix_retained_system} and the modal representation \eqref{eq:appendix_mild_solution}.

For every $t>0$, the homogeneous component satisfies $\sum_{n=1}^{\infty}e^{-2(\lambda_n-q_c)t}|z_n^0|^2<\infty$.

Furthermore, by \eqref{eq:appendix_input_bound},
\begin{equation*}
    \sum_{n=N}^{\infty}
    \left|
        \beta_n\int_0^t
        e^{-(\lambda_n-q_c)(t-s)}u(s)\,ds
    \right|^2
    \leq
    4M_\pi^2
    \sum_{n=N}^{\infty}\frac{\beta_n^2}{\lambda_n^2}
    <\infty.
\end{equation*}

Therefore, $z(t)\in L^2(0,1)$ for every $t>0$.

Finally, the global Lipschitz continuity of $\pi$, together with the preceding estimates, gives continuous dependence on the initial condition. This concludes the proof.
\hfill$\square$
\end{pf}

It remains to relate the modal solution to the original boundary-controlled PDE. Since $Z_{n_0}\in C^{1}([0,\infty);\mathbb R^{n_0})$ and $\pi$ is globally Lipschitz, the feedback input $u=\pi(Z_{n_0})$ belongs to $W_{\mathrm{loc}}^{1,\infty}(0,\infty)$.

Let $T(t)$ denote the semigroup generated by $-\mathcal A+q_cI$ on $\mathcal H_{-1/2}$. With the lifting introduced in \eqref{eq:lift}, set $w^0=z^0+b\,u(0)$ and define
\begin{equation*}
    w(t)=T(t)w^0+ \int_0^t T(t-s) \bigl(a\,u(s)+b\,\dot u(s)\bigr)\,ds.
\end{equation*}

Then $w\in C([0,\infty);\mathcal H_{-1/2})$ is the unique mild solution of the homogeneous-boundary lifted system.

Consequently, $z(t)=w(t)-b\,u(t)$ defines the unique lifted mild solution of \eqref{eq:pde_main}--\eqref{eq:pde_ini}.

Projecting this identity onto $\Phi_n$ and integrating the term containing $\dot u$ by parts yields \eqref{eq:appendix_closed_loop}, because $\beta_n=a_n+(-\lambda_n+q_c)b_n$. Hence the modal solution constructed above is precisely the lifted mild solution of the original boundary-controlled PDE.

\bibliographystyle{cas-model2-names}

\bibliography{cas-refs}

@book{curtain2020,
	author    = {Curtain, Ruth F. and Zwart, Hans},
	title     = {Introduction to Infinite-Dimensional Systems Theory: A State-Space Approach},
	series    = {Texts in Applied Mathematics},
	volume    = {71},
	publisher = {Springer},
	address   = {Cham},
	year      = {2020},
	edition   = {2}
}

@book{bastin2016,
	title     = {Stability and boundary stabilization of 1-{D} hyperbolic systems},
	author    = {Bastin, G. and Coron, J.-M.},
	year      = {2016},
	publisher = {Birkh\"{a}user}
}

@book{krstic2008,
	title     = {Boundary control of {PDE}s: {A} course on backstepping designs},
	author    = {Krstic, M. and Smyshlyaev, A.},
	year      = {2008},
	publisher = {SIAM}
}

@article{yin2021imitation,
	title   	= {Imitation learning with stability and safety guarantees},
	author  	= {Yin, He and Seiler, Peter and Jin, Ming and Arcak, Murat},
	journal 	= {IEEE Control Systems Letters},
	volume  	= {6},
	pages   	= {409--414},
	year    	= {2021},
	publisher	= {IEEE}
}

@article{lhachemi2020pi,
	title   	= {{PI} regulation of a reaction--diffusion equation with delayed boundary control},
	author  	= {Lhachemi, Hugo and Prieur, Christophe and Tr{\'e}lat, Emmanuel},
	journal	 	= {IEEE Transactions on Automatic Control},
	volume  	= {66},
	number  	= {4},
	pages   	= {1573--1587},
	year    	= {2020},
	publisher = {IEEE}
}

@inproceedings{de2023imitation,
	title   		= {Imitation learning control for thermoacoustic stabilization of a {R}ijke tube},
	author  		= {de Andrade, Gustavo Artur and Fiacchini, Mirko and Prieur, Christophe},
	booktitle		= {2023 62nd IEEE Conference on Decision and Control (CDC)},
	pages   		= {2958--2963},
	year    		= {2023},
	organization= {IEEE}
}

@book{boyce2021elementary,
	title   	= {Elementary differential equations and boundary value problems},
	author  	= {Boyce, William E and DiPrima, Richard C and Meade, Douglas B},
	year    	= {2021},
	publisher	= {John Wiley \& Sons}
}

@book{tucsnak2009observation,
	title     = {Observation and control for operator semigroups},
	author    = {Tucsnak, Marius and Weiss, George},
	year      = {2009},
	publisher = {Springer}
}

@article{luo2025physics,
	author  = {Luo, Kuang and Zhao, Jingshang and Wang, Yingping and Li, Jiayao and Wen, Junjie and Liang, Jiong and Soekmadji, Henry and Liao, Shaolin},
	journal = {Artificial Intelligence Review},
	number  = {10},
	pages   = {1--43},
	title   = {Physics-informed neural networks for {PDE} problems: a comprehensive review},
	volume  = {58},
	year    = {2025}
}

@article{mattioni2023enhancing,
	author  = {Mattioni, Andrea and Zoboli, Samuele and Mavkov, Bojan and Astolfi, Daniele and Andrieu, Vincent and Witrant, Emmanuel and Frasca, Paolo and Prieur, Christophe},
	journal = {Fusion Engineering and Design},
	pages   = {114008},
	publisher = {Elsevier},
	title   = {Enhancing deep reinforcement learning with integral action to control tokamak safety factor},
	volume  = {196},
	year    = {2023}
}

@article{tanyu2023deep,
	author  = {Tanyu, Derick Nganyu and Ning, Jianfeng and Freudenberg, Tom and Heilenk{\"o}tter, Nick and Rademacher, Andreas and Iben, Uwe and Maass, Peter},
	journal = {Inverse Problems},
	number  = {10},
	pages   = {103001},
	title   = {Deep learning methods for partial differential equations and related parameter identification problems},
	volume  = {39},
	year    = {2023}
}

@article{zhang2026operator,
  title={Operator learning for robust stabilization of linear {M}arkov-jumping hyperbolic {PDE}s},
  author={Zhang, Yihuai and Auriol, Jean and Yu, Huan},
  journal={Automatica},
  volume={185},
  pages={112809},
  year={2026},
  publisher={Elsevier}
}

@book{morris2020control,
	author    = {Kirsten A. Morris},
	title     = {Control of Systems Governed by Partial Differential Equations},
	year      = {2020},
	publisher = {Springer},
	address   = {Cham, Switzerland},
	series    = {Applied Mathematical Sciences},
	volume    = {211},
}

@book{li1995optimal,
	author    = {Xunjing Li and Jiongmin Yong},
	title     = {Optimal Control Theory for Infinite Dimensional Systems},
	year      = {1995},
	publisher = {Birkhäuser Boston},
	address   = {Boston, MA},
	series    = {Systems \& Control: Foundations \& Applications},
}

@article{pinguet2023data,
	title   	= {A data-based neural controller training method with tunable stability margin using multi-objective optimization},
	author  	= {Pinguet, J{\'e}r{\'e}my and Feyel, Philippe and Sandou, Guillaume},
	journal 	= {IFAC-PapersOnLine},
	volume  	= {56},
	number  	= {2},
	pages   	= {3092--3099},
	year    	= {2023},
	publisher	= {Elsevier}
}

@article{zare2024survey,
	title   	= {A survey of imitation learning: Algorithms, recent developments, and challenges},
	author  	= {Zare, Maryam and Kebria, Parham M and Khosravi, Abbas and Nahavandi, Saeid},
	journal 	= {IEEE Transactions on Cybernetics},
    volume  	= {54},
	number  	= {12},
	pages   	= {7173--7186},
	year    	= {2024},
}

@inproceedings{bhan2023operator,
	author    = {Bhan, Luke and Shi, Yuanyuan and Krstic, Miroslav},
	title     = {Operator Learning for Nonlinear Adaptive Control},
	booktitle = {Proceedings of the 5th Annual Learning for Dynamics and Control Conference},
	series    = {Proceedings of Machine Learning Research},
	volume    = {211},
	pages     = {346--357},
	year      = {2023}
}

@article{bhan2024neuraloperators,
	author  = {Bhan, Luke and Shi, Yuanyuan and Krstic, Miroslav},
	title   = {Neural Operators for Bypassing Gain and Control Computations in {PDE} Backstepping},
	journal = {IEEE Transactions on Automatic Control},
	volume  = {69},
	number  = {8},
	pages   = {5310--5325},
	year    = {2024},
}

@article{krstic2024neuraloperators,
	author  = {Krstic, Miroslav and Bhan, Luke and Shi, Yuanyuan},
	title   = {Neural Operators of Backstepping Controller and Observer Gain Functions for Reaction--Diffusion {PDEs}},
	journal = {Automatica},
	volume  = {164},
	pages   = {111649},
	year    = {2024},
}

@article{bhan2025adaptive,
	author  = {Bhan, Luke and Shi, Yuanyuan and Krstic, Miroslav},
	title   = {Adaptive Control of Reaction--Diffusion {PDEs} via Neural Operator-Approximated Gain Kernels},
	journal = {Systems \& Control Letters},
	volume  = {195},
	pages   = {105968},
	year    = {2025},
}

@inproceedings{bhan2024pdecontrolgym,
	author    = {Bhan, Luke and Bian, Yuexin and Krstic, Miroslav and Shi, Yuanyuan},
	title     = {{PDE} Control Gym: A Benchmark for Data-Driven Boundary Control of Partial Differential Equations},
	booktitle = {Proceedings of the 6th Annual Learning for Dynamics and Control Conference},
	series    = {Proceedings of Machine Learning Research},
	volume    = {242},
	pages     = {1083--1095},
	year      = {2024}
}

@article{lu2021deeponet,
	author  = {Lu, Lu and Jin, Pengzhan and Pang, Guofei and Zhang, Zhongqiang and Karniadakis, George Em},
	title   = {Learning Nonlinear Operators via {DeepONet} Based on the Universal Approximation Theorem of Operators},
	journal = {Nature Machine Intelligence},
	volume  = {3},
	pages   = {218--229},
	year    = {2021},
}

@inproceedings{li2021fourier,
	author    = {Li, Zongyi and Kovachki, Nikola and Azizzadenesheli, Kamyar and Liu, Burigede and Bhattacharya, Kaushik and Stuart, Andrew and Anandkumar, Anima},
	title     = {Fourier Neural Operator for Parametric Partial Differential Equations},
	booktitle = {International Conference on Learning Representations},
	year      = {2021}
}

@article{kovachki2023neural,
	author  = {Kovachki, Nikola and Li, Zongyi and Liu, Burigede and Azizzadenesheli, Kamyar and Bhattacharya, Kaushik and Stuart, Andrew and Anandkumar, Anima},
	title   = {Neural Operator: Learning Maps Between Function Spaces with Applications to {PDEs}},
	journal = {Journal of Machine Learning Research},
	volume  = {24},
	number  = {89},
	pages   = {1--97},
	year    = {2023}
}

@article{east2022imitation,
	author  = {East, Sebastian},
	title   = {Imitation Learning of Stabilizing Policies for Nonlinear Systems},
	journal = {European Journal of Control},
	volume  = {68},
	pages   = {100678},
	year    = {2022},
}

@article{yin2022stability,
	author  = {Yin, He and Seiler, Peter and Arcak, Murat},
	title   = {Stability Analysis Using Quadratic Constraints for Systems with Neural Network Controllers},
	journal = {IEEE Transactions on Automatic Control},
	volume  = {67},
	number  = {4},
	pages   = {1980--1987},
	year    = {2022},
}

@article{Lhachemi26,
	author  = {Lhachemi, Hugo and Munteanu, Ionuţ and Prieur, Christophe},
	journal = {IEEE Transactions on Automatic Control}, 
	title   = {Boundary Output Feedback Stabilization of Nonlinear Diffusion Equations}, 
	year    = {2026},
	volume  = {71},
	number  = {7},
	pages   = {4764-4771},
}

\end{document}